\documentclass[12pt,a4paper]{article}
\usepackage{a4wide,amsmath,amsthm,epsfig,graphicx,psfrag}
\usepackage{amsmath,amsthm}
\usepackage{amsfonts}
\usepackage{amssymb}
\usepackage{mathscinet}
\usepackage{enumerate}
\usepackage[latin1]{inputenc}  

\newcommand{\url}[1]{}

\newcommand{\Px}{\Bbb{P}}
\newcommand{\E}{\Bbb{E}}

\newcommand{\R}{{\mathbb R}}
\newcommand{\Ex}{\Bbb{E}}

\newcommand{\Hs}{\mathcal{H}_{\text{\bf vol}}}

\newcommand{\tsigma}{\tilde{\sigma}}

\newcommand{\xsd}{\xi^*}
\newcommand{\ysd}{y^*}
\newcommand{\xss}{x^*_\sigma}
\newcommand{\yse}{\ysd_\eta}
\newcommand{\ysp}{\Upsilon^*}

\newcommand{\uteta}{\utilde{\eta}}

\newcommand{\uh}[1]{ \underset{ \check{}} {#1 }}

\newtheorem{theorem}{Theorem}[section]
\newtheorem{proposition}[theorem]{Proposition}
\newtheorem{lemma}[theorem]{Lemma}
\newtheorem{corollary}[theorem]{Corollary}

\newtheorem{remark}[theorem]{Remark}

\newtheorem{convention}[theorem]{Convention}

\def\title#1{\hfil\break\hfil\break
\hfil\break\par\addvspace\baselineskip\noindent
\ignorespaces{\LARGE\bf#1}\hfil\break}

\def\author#1{\par\addvspace\baselineskip\noindent
\ignorespaces{\large\bf#1}}

\def\institute#1{\par\addvspace\baselineskip\noindent
\ignorespaces{\small#1}\hfil\break}

\begin{document}

\title{{\bf A Call-Put Duality for Perpetual American Options }}
\author{Aur\'{e}lien Alfonsi and Benjamin Jourdain}
\institute{  CERMICS, projet MATHFI, Ecole Nationale des Ponts et Chauss\'{e}es, 6-8
avenue Blaise Pascal, Cit\'{e} Descartes, Champs sur Marne,
77455 Marne-la-vall\'{e}e, France.\\ e-mail : {\tt \{alfonsi,jourdain\}@cermics.enpc.fr } \\
\today }

{\abstract It is well known \cite{Dup}, \cite{AC} that in models with time-homogeneous local
volatility functions and constant interest and dividend
rates, the European Put prices are transformed into European Call prices
by the simultaneous exchanges of
the interest and dividend rates and of the strike and spot price of the underlying.
This paper investigates such a Call Put duality for perpetual
American options. It turns out that the perpetual American Put price
is equal to the perpetual American Call price in a model
where, in addition to the previous exchanges between the spot price and the
strike and between the interest and dividend rates, the local volatility
function is modified. We prove that equality of the dual volatility functions only holds in
the standard Black-Scholes model with constant volatility. Thanks to these duality
results, we design a theoretical calibration procedure of
the local volatility function from the perpetual Call
and Put prices for a fixed spot price
$x_0$. The knowledge of the Put (resp. Call) prices for all strikes
enables to recover the local
volatility function on the interval $(0,x_0)$ (resp. $(x_0,+\infty)$). }

\vspace{0.2cm}
\thispagestyle{empty}

\section*{Introduction}
In a model with local volatility function $\varsigma(t,x)$, interest rate $r$ and dividend rate $\delta$
\begin{equation}
   \begin{cases}
      dS^x_t=\varsigma(t,S^x_t)S^x_tdW_t+(r-\delta)S^x_tdt,\;t\geq 0\\
S^x_0=x
   \end{cases}
\label{intedsprim}\end{equation} the initial price
$$h(T,y)=\E\left[e^{-rT}(y-S^x_T)^+\right]$$
of the European Put option considered as a function of the maturity $T>0$ and the Strike $y>0$ solves Dupire's partial differential equation \cite{Dup} :
\begin{equation*}
   \begin{cases}
      \partial_T
h(T,y)=\frac{\varsigma^2(T,y)y^2}{2}\partial^2_{yy}h(T,y)+(\delta-r)y\partial_{y}h(T,y)-\delta
h(T,y),\;T,y>0\\
h(0,y)=(y-x)^+,\;y>0
   \end{cases}
\end{equation*}
One easily deduces that the function $h(T-t,y)$ for $(t,y)\in [0,T] \times \R_+^*$
satisfies the pricing partial differential equation for the Call
option with strike $x$ and maturity $T$ in the model
\begin{equation}
   \begin{cases}
      d\bar{S}^{y,T}_t=\varsigma(T-t,\bar{S}^{y,T}_t)\bar{S}^{y,T}_tdW_t+(\delta-r)\bar{S}^{y,T}_tdt,\;t\in[0,T]\\
\bar{S}^{y,T}_0=y
   \end{cases}
\end{equation}
with local volatility function $\varsigma(T-t,y)$, interest rate $\delta$ and
dividend rate $r$. Therefore $h(T,y)=\E\left[e^{-\delta
    T}(\bar{S}^{y,T}_T-x)^+\right]$ and one deduces the following Call-Put duality
relation which is also a consequence of \cite{AC}
$$\forall T\geq 0,\;\forall x,y>0,\;\E\left[e^{-rT}(y-S^x_T)^+\right]=\E\left[e^{-\delta
    T}(\bar{S}^{y,T}_T-x)^+\right].$$
Since it derives from Dupire's formula, this Call-Put duality equality
    is closely related to calibration issues. One remarks that in the particular case of a time-homogeneous volatility function
    ($\varsigma(t,x)=\sigma(x)$), then $\bar{S}^{y,T}_t$ also evolves
    according to the same time-homogeneous volatility function.\par
In this work, we are interested in deriving such a Call-Put duality
    relation in the case of American options and in investigating
    consequences in terms of calibration. In the Black-Scholes model with
    constant volatility $\varsigma(t,x)=\sigma$, when $\tau$ denotes a
    bounded stopping-time of the natural filtration of the Brownian
    motion $(W_t)_{t\geq 0}$, one has
\begin{align*}\E\left[e^{-r\tau}\left(y-xe^{\sigma
    W_\tau+(r-\delta-\frac{\sigma^2}{2})\tau}\right)^+\right]&=\E\left[e^{-\delta
    \tau}e^{\sigma W_\tau-\frac{\sigma^2}{2}\tau}\left(ye^{-\sigma
    W_\tau+(\delta-r+\frac{\sigma^2}{2})\tau}-x\right)^+\right]
\\
&=\E\left[e^{-\delta
    \tau}\left(ye^{-\sigma
    W_\tau+(\delta-r-\frac{\sigma^2}{2})\tau}-x\right)^+\right]
\end{align*}
where the second equality follows from Girsanov theorem. Taking the
supremum over all stopping-times $\tau$ smaller than $T$ one deduces
that the price of the American Put option with maturity $T$ is equal to
the price of the American Call option with the same maturity up to the
simultaneous exchange between the underlying spot price and the strike
and between the interest and dividend rates. 
Extensions of this result when the underlying evolves according to the
exponential of a L\'evy process have been obtained in
\cite{FM}. Let us also mention that another kind of duality has been investigated
in~\cite{PS}. But, to our knowledge, no study has been
devoted to the case of models with local volatility functions like
\eqref{intedsprim}.\par
In the present paper, we consider the case of
perpetual ($T=+\infty$) American options in models with time-homogeneous local
volatility functions $\varsigma(t,x)=\sigma(x)$.
In the first part, we recover well-known properties of the perpetual American call
and put pricing functions by extending an approach recently developed by Beibel and
Lerche \cite{BL} in the Black-Scholes case. This makes the paper self-contained.\\ In the second part, we introduce the
framework used in the remaining of the paper. \\
In the third part of the paper, we 
consider the exercise boundaries 
as functions of the strike variable and characterize them as the unique
solutions of some non-autonomous ordinary
differential equations. \\
The fourth part is dedicated to our main result. We prove that the
perpetual American Put prices are equal to the perpetual American Call
prices in a model where, in addition to the exchanges between the spot
price of the underlying and the strike and between the interest and dividend
rates, the volatility function is modified. We also derive an expression
of this modified volatility function. Notice that in the European
case presented above, time-homogeneous volatility functions are not modified.
\\
The fifth part addresses calibration issues. It turns out that for a
given initial value $x_0>0$ of the underlying one recovers the
restriction of the time-homogeneous volatility function
$\sigma(x)$ to $(0,x_0]$ (resp. $[x_0,+\infty)$) from the perpetual Put
(resp. Call)
prices for all strikes.\\
In the last part, we show that at least when $\delta<r$, in the class of
volatility functions analytic in a neighbourhood of the origin, the
only ones invariant by our duality result are the constants. This means
that the case of the standard Black-Scholes model presented above is very specific.

{\bf Acknowledgements. } We thank Damien Lamberton (Univ. Marne-la-vallée) and Mihail
Zervos (King's College) for interesting
discussions. We also thank Alexander Schied (TU Berlin) for pointing out the work  of
Beibel and Lerche \cite{BL} to us and Antonino Zanette (University of Udine) for providing us with
the routine that calculates American option prices.

\section{Perpetual American put and call pricing}\label{Sec_FB}

We consider a constant interest spot-rate $r$ that is assumed to
be nonnegative and an asset $S_t$ which pays a constant dividend rate $\delta \ge 0$ and is
driven by a homogeneous volatility function $\sigma: \mathbb{R}^*_+ \rightarrow
\mathbb{R}^*_+$ that satisfies the following hypothesis.  

{\bf Hypothesis $(\Hs)$:} {\em $\sigma$ is continuous on $\mathbb{R}^*_+$ and there are  $0<\underline{\sigma}<\overline{\sigma}<+\infty$ such that:
  \[ \forall x>0, \underline{\sigma} \le \sigma(x) \le
  \overline{\sigma}. \]
}
In other words, $S_t$ is assumed to follow under the risk-neutral measure the SDE:
\begin{equation}\label{EDS_stock}
  dS_t=S_t((r-\delta)dt+\sigma(S_t)dW_t).
\end{equation}
With the assumption made on $\sigma$, we
know that for any initial condition $x \in \mathbb{R}^*_+$, there is a unique solution in the sense of probability law (see for example
Theorem 5.15 in \cite{KS}, using a log transformation) denoted by $(S^{x}_t,t \ge 0)$.
Moreover, Theorem 4.20 ensures that the strong Markov property holds for $(S^x_t,t \ge 0)$. Under that model, we
denote by
\[ P_\sigma (x,y)=\underset{\tau \in
  \mathcal{T}_{0,\infty}} {\sup} \Ex \left[e^{-r\tau}(y-S^{x}_{\tau})^+
\right] \text{ and } C_\sigma (x,y)= \underset{\tau \in
  \mathcal{T}_{0,\infty}} {\sup} \Ex \left[e^{-r\tau}(S^{x}_{\tau}-y)^+
\right]  \]
respectively the prices of the American perpetual put and call options with strike $y>0$ and
spot~$x$.
Here, $\mathcal{T}_{0,\infty}$ simply denotes the set of the stopping times
with respect to the natural filtration of $(S^{x}_t,t \ge 0)$. Since $e^{-rt}S^{x}_t=x \exp \left(-\delta t +
\int_0^t \sigma (S^{x}_u) dW_u - \frac{1}{2}\int_0^t \sigma^2 (S^{x}_u)du \right)$
and $\int_0^t \sigma^2 (S^{x}_u)du \ge \underline{\sigma}^2 t
\underset{t \rightarrow + \infty}{\rightarrow} + \infty$, it follows from the
Dubins-Schwarz theorem that 
\begin{equation}\label{Dub}
  e^{-rt}S^{x}_t \underset{t \rightarrow + \infty}{\rightarrow} 0 \text{
    a.s.}\end{equation}
As a consequence, $e^{-r  t}(y-S^{x}_{ t})^+\underset{t \rightarrow + \infty}
{\rightarrow}\mathbf{1}_{\{r=0\}}y$ and $e^{-r  t}(S^{x}_{ t}-y)^+\underset{t \rightarrow + \infty}
{\rightarrow}0$ a.s. On $\{ \tau= \infty \}$, we thus set
\begin{equation}\label{tau_infty}
e^{-r\tau}(y-S^{x}_{\tau})^+=\mathbf{1}_{\{r=0\}}y \text{ and }
e^{-r\tau}(S^{x}_{\tau}-y)^+=0.
\end{equation}

Let us  consider the second-order ordinary differential equation
 \begin{equation}\label{EDf}
   \frac{1}{2} \sigma^2(x) x^2 f''(x)+ (r-\delta) x  f'(x) -  rf(x) =0, \ x>0.
 \end{equation}
According to Borodin and Salminen 
(\cite{BS}, chap. 2) the functions \begin{equation} \label{def_f_ud}
\forall x>0, \   f_\uparrow(x)=\left\{\begin{array}{l} \Ex[e^{-r\tau^x_1}] , \text{ if } x \le 1 \\
  1/\Ex[e^{-r\tau^1_x}], \text{ if } x > 1
\end{array}\right. \text{ and }f_\downarrow(x)=\left\{\begin{array}{l}  1/\Ex[e^{-r\tau^1_x}] , \text{ if } x \le 1 \\
\Ex[e^{-r\tau^x_1}] , \text{ if } x > 1
\end{array}\right. 
\end{equation}
where for $x,y>0$, $ \tau^x_y=\inf\{t \ge 0, S^x_t=y \} \ \  (\ \inf \emptyset=+\infty \ )$,
are the unique solutions (up to a multiplicative
constant) that are positive and respectively increasing and 
decreasing. The volatility function $\sigma$ being continuous, (\ref{def_f_ud})
ensures that these functions are $\mathcal{C}^2$ on $\mathbb{R}^*_+$.

\begin{remark}\label{Tps_atteinte} It is easy using the strong Markov property to
  get:
  \[ \forall x,y>0, \  \Ex [e^{-r \tau^x_y}] = \begin{cases}
    f_\uparrow(x)/ f_\uparrow(y),\text{ if } x \le y\\
    f_\downarrow(x)/f_\downarrow(y) ,\text{ if } x \ge y.
  \end{cases}\]

\end{remark}
\begin{remark}\label{LIMCOND}
Assuming $r>0$,  one has (see Borodin and Salminen \cite{BS})
  \begin{eqnarray*}
&&  \underset{x \rightarrow 0}  {\lim} f_\downarrow(x)=+\infty, \  
  \underset{x \rightarrow +\infty }  {\lim} f_\downarrow(x)=   0\\
&& \underset{x \rightarrow 0}  {\lim} f_\uparrow(x)= 0 ,\ 
  \underset{x \rightarrow +\infty }  {\lim} f_\uparrow(x)= +\infty.\\
\end{eqnarray*}
  The function $f_\downarrow$ (resp. $f_\uparrow$) is thus, up to a multiplicative constant,
  the unique solution to~(\ref{EDf})  such that $\underset{x \rightarrow +\infty}  {\lim} f(x)=0$ (resp. $\underset{x
    \rightarrow 0}  {\lim} f(x)=0$).
\end{remark}
\begin{remark}\label{solu_delta0}
In the case $\delta=0$, we  have the analytical solutions~:
\[ f_\downarrow(x)=\frac{\varphi(x)}{\varphi(1)}\;\mbox{where}\;\varphi(x)=x \int^{+\infty}_x\left( \frac{1}{v^2} \exp \left[- \int_1^v
      \frac{2r}{u\sigma^2(u)}du\right] \right) dv,\  
  f_ \uparrow(x)=x. \]
 Indeed,  since $f(x)=x$ is
solution of $\frac{1}{2} \sigma^2(x) x^2f''(x)+rxf'(x)-rf(x)=0$, we search a general solution that can be written
$f(x)=x \tilde{f}(x)$. This leads to $ \frac{1}{2} \sigma^2(x) x
\tilde{f}''(x)+ (r  +\sigma^2(x)  )  \tilde{f}'(x)  =0$ and then
$\tilde{f}'(x)=\frac{C_1}{x^2} \exp \left[- \int_1^x
  \frac{2r}{u\sigma^2(u)}du \right]$. Therefore, $
\exists C_1,C_2
\in \mathbb{R}, \ f(x)=C_2 x + C_1 x \int^{+\infty}_x\left( \frac{1}{v^2} \exp \left[- \int_1^v
  \frac{2r}{u\sigma^2(u)}du\right] \right) dv. $
\end{remark}

Now, we are in position to show the existence of an optimal stopping time and give
 the call and put prices. Let us mention here that the problem of perpetual optimal
stopping is treated in the paper of Dayanik and Karatzas
\cite{DK}, for a general payoff function and an underlying evolving according to a
general one-dimensional time homogeneous diffusion process.
Villeneuve~\cite{Villeneuve} considers a model where the constant dividend rate
$\delta$ in~(\ref{EDS_stock}) is replaced by a function $\delta(S_t)$ and gives
sufficient condition on the payoff function ensuring that a threshold strategy is optimal.
Here, we give a direct proof that generalizes the approach developed by Beibel
and Lerche~\cite{BL} in the Black-Scholes case. 
\begin{theorem}\label{Opt_stop_put}
  Assume $r>0$. For any strike $y>0$, there is a unique $\xss(y) \in (0,y) \cap (0,\frac{r}{\delta}y]$
 such that $\tau^P_x=\inf \{t \ge 0, S^x_t \le \xss(y) \}$ (convention $\inf
 \emptyset=+\infty$) is an optimal stopping time for the put and:
 \begin{equation}\label{DEFx*}
  \forall x \le \xss(y), P_\sigma(x,y)=(y-x)^+, \ \forall x >
  \xss(y), P_\sigma(x,y)= \frac{y- \xss(y) }{
    f_\downarrow (\xss (y))}f_\downarrow (x) > (y-x)^+.  \hspace{0.1cm}
\end{equation}
In addition,  we have $ f'_\downarrow(\xss(y) )<0 $ and:
 \begin{equation}\label{DEF2x*}
 \xss(y)-y=\frac{f_\downarrow (\xss (y))}{ f'_\downarrow(\xss(y) )}.
\end{equation}
Last, the smooth-fit principle holds: $\partial_x P_\sigma(x^*_\sigma(y),y)=-1$.
\end{theorem}
\begin{remark} If $r=0$, $\forall x,y>0, P_\sigma(x,y)=y$ since for any stopping time
  $\tau$, $(y-S^x_{\tau})^+\le y$ and equality holds for $\tau=+\infty$ by~(\ref{tau_infty}).
\end{remark}


\begin{proof}
  Let us define:
  \[ \forall z>0,\  h(z)=\frac{(y-z)^+}{f_\downarrow (z)} \text{ and }
  h^*=\underset{z>0}{\sup} \ h(z).\]
  Since the function $h$ is continuous such that $h(y)=0$ and $h(0+)=0$
  (Remark~\ref{LIMCOND}), $\xss(y) = \sup \{z>0 , \ h(z)=h^*\}$ belongs to $(0,y)$
  and is such that $h(\xss(y))=h^*$.
  Since the function $h$ is
  $\mathcal{C}^2$ on $(0,y)$, we have $h'(\xss(y))=0$ and $h''(\xss(y)) \le 0$. These
  conditions give easily\[ f_\downarrow (\xss (y))+ (y-\xss(y) )
  f'_\downarrow(\xss(y) )=0 \text{ and } f''_\downarrow(\xss(y) ) \ge 0.
  \]
  Since $f_\downarrow$ is positive and $\xss(y)<y$, we have  $f'_\downarrow(\xss(y)
  )<0$ and deduce~(\ref{DEF2x*}). The second order condition and equation~(\ref{EDf})
  then give $\xss(y)
  (r-\delta)f'_\downarrow(\xss(y)  ) -rf_\downarrow(\xss(y) )  \le 0$ and
  so \[ry - \delta \xss(y) \ge 0.\]

  Now let us check the optimality of $\tau^P_x$ and  consider $\tau \in
  \mathcal{T}_{0,\infty}$. By Fatou's lemma and Doob's optional sampling theorem,
  we have
  \begin{eqnarray*}
    \Ex[e^{-r \tau}(y-S^{x}_{\tau})^+ ] &\le& \underset{t \rightarrow
    +\infty}{\liminf}\Ex \left[ e^{-r\tau \wedge t} (y- S^{x}_{\tau \wedge
      t})^+   \right]\\
  &=&\underset{t \rightarrow
    +\infty}{\liminf}\Ex \left[ e^{-r\tau \wedge t}  f_\downarrow(S^{x}_{\tau \wedge
      t}) h(S^{x}_{\tau \wedge
      t})   \right] \\
  &\le &h(\xss(y))\underset{t \rightarrow
    +\infty}{\liminf}\Ex \left[ e^{-r\tau \wedge t}  f_\downarrow(S^{x}_{\tau \wedge
      t})   \right] \le h(\xss(y))f_\downarrow(x)
\end{eqnarray*}
since $e^{-rt}f_\downarrow(S^{x}_t)=f_\downarrow(x)+\int_0^t
  e^{-ru}\sigma(S^{x}_u)S^{x}_u f_\downarrow '(S^{x}_u) dW_u$ is a
  nonnegative local martingale and therefore a supermartingale.  If $x \ge \xss(y)$,
  we have using Remark \ref{Tps_atteinte}:
\[\Ex[e^{-r \tau^P_x}(y-S^{x}_{\tau^P_x})^+ ]= \Ex[e^{-r
  \tau^x_{\xss (y)}}(y-S^{x}_{\tau^x_{\xss(y)}})^+ ]=(y-\xss (y))\Ex[e^{-r \tau^x_{\xss
    (y)}}]=h(\xss(y))f_\downarrow(x)\]
and $\tau^P_x$ is optimal for $x \ge \xss(y)$. Since  $\xss(y) = \sup \{z>0 , \
h(z)=h^*\}$, we have $(y-x)^+=h(x)f_\downarrow(x)<f_\downarrow(x)h(\xss(y))$ for $x > \xss(y)$, and finally deduces (\ref{DEFx*}) for $x \ge \xss(y)$.

We consider now the complementary case $x \in (0,\xss(y))$, and set $\tau \in
\mathcal{T}_{0,\infty}$. Using the strong Markov property and the optimality result
when the initial spot is $\xss(y)$, we get
\[\Ex[e^{-r \tau}
  (y-S^{x}_{\tau})^+] \le \Ex[e^{-r \tau \wedge \tau^x_{\xss (y)}}
  (y-S^{x}_{\tau \wedge \tau^x_{\xss (y)} })^+]. \]
On $\{t < \tau^x_{\xss (y)} \}$, we have $S^{x}_t<\xss(y)$,
$de^{-rt}(y-S^{x}_t)=e^{-rt}\underset{\le 0}{\underbrace{(\delta S^{x}_t -ry )}}dt
-e^{-rt}\sigma(S^{x}_t)S^{x}_t dW_t$ and so $\Ex[e^{-r \tau \wedge \tau^x_{\xss (y)}}
  (y-S^{x}_{\tau \wedge \tau^x_{\xss (y)} })^+]\le \underset{t \rightarrow
    +\infty}{\liminf} \Ex[e^{-r \tau \wedge \tau^x_{\xss (y)} \wedge t}
  (y-S^{x}_{\tau \wedge \tau^x_{\xss (y)} \wedge t})^+] \le (y-x)$. 
\end{proof}

Now, we state the similar result for the call prices.

\begin{theorem}\label{Opt_stop_call}
  Assume $\delta>0$. For any strike $y>0$, there is a unique $\ysp_{\sigma}(y) \in (y,\infty)
  \cap [\frac{r}{\delta}y,+\infty)$ such that  $\tau^C_x=\inf \{t \ge 0, S^x_t \ge
  \ysp_{\sigma}(y) \}$ is an optimal stopping time for the call and:
 \begin{equation} \label{DEFy*} 
  \forall x \ge \ysp_{\sigma}(y), C_\sigma(x,y)=(x-y)^+, \   \forall x <
  \ysp_{\sigma}(y), C_\sigma(x,y)=\frac{\ysp_{\sigma}(y)-y }{
    f_\uparrow (\ysp_{\sigma} (y))}f_\uparrow (x) > (x-y)^+.
\end{equation}
In addition, we have $f'_\uparrow(\ysp_{\sigma}(y) )>0$ and:
 \begin{equation}\label{DEF2y*}
\ysp_{\sigma}(y)-y=\frac{f_\uparrow (\ysp_{\sigma} (y))}{ f'_\uparrow(\ysp_{\sigma}(y) )}.
\end{equation}
Last, the smooth-fit principle holds: $\partial_x C_\sigma(\ysp(y),y)=1$.
\end{theorem}
\begin{remark} If $\delta=0$, $\forall x,y>0, C_\sigma(x,y)=x$. Indeed, the
Call-Put parity $\Ex[e^{-rt}(S^{x}_t-y)^+]=x-ye^{-rt}+\Ex[e^{-rt}(y-S^{x}_t)^+]$
gives the convergence to $x$ in both cases $r>0$ and $r=0$ when $t \rightarrow +\infty$.
Now, thanks to the Fatou lemma, we have for $\tau \in \mathcal{T}_{0,\infty}$:
$  \Ex \left[e^{-r\tau}(S^{x}_{\tau}-y)^+ \right]  \le  \underset{t \rightarrow
    +\infty}{\liminf}\Ex \left[ e^{-r\tau \wedge t}(S^{x}_{\tau \wedge t}-y)^+
  \right] 
  \le  \underset{t \rightarrow
    +\infty}{\liminf}\Ex \left[ e^{-r\tau \wedge t}S^{x}_{\tau \wedge t}
  \right]=x. $ 
\end{remark}
\begin{proof}
The proof works as for the put, and we just hint the differences. We define
  \[ \forall z>0, h(z)=\frac{(z-y)^+}{f_\uparrow (z)} \text{ and }
  h^*=\underset{z>0}{\sup} \ h(z).\]
 Let us admit for a while that $\underset{z\rightarrow + \infty}{\lim }   h(z)=0$.
 Then, since $h(y)=0$ and $h$ is continuous, $h$ reaches its maximum in $\ysp_{\sigma}(y)=\inf \{z>0, \
 h(z)=h^*\}$, and $\ysp_{\sigma}(y) \in (y,\infty)$. This gives (\ref{DEF2y*}). We then consider the case $x \le \ysp_{\sigma}(y)$
 and show that $\tau^C_x=\tau^x_{\ysp_{\sigma}(y)}$ is optimal. Note that in the special case
 $r=0$, we have to use Proposition \ref{solu0_fg} which is stated below. Finally, we prove that $\tau^C_x$ is
 optimal when $x > \ysp_{\sigma}(y)$ using that $\delta\ysp_{\sigma}(y) -r y \ge 0 $.

 Now, let us check that $\underset{z\rightarrow + \infty}{\lim }   h(z)=0$. In the
 case $r=0$, it is straightforward using the explicit form given in Proposition
 \ref{solu0_fg} below that $f_\uparrow(x) 
 \ge\frac{1}{1+\frac{ 2 \delta }{ \overline{\sigma}^2}} x^{1+ 2 \delta /
   \overline{\sigma}^2}-1$ for $x \ge 1$, and we
   have then $\underset{z\rightarrow + \infty}{\lim }   h(z)=0$. When $r>0$, Itô's Formula gives
\[de^{-r t}(S^{1}_t)^{1+a}=e^{-r t}(S^{1}_t)^{1+a}\left\{(a+1)\sigma(S^{1}_t)dW_t
  +[a(r+(a+1)\sigma^2(S^{1}_t)/2 ) - (a+1) \delta ]dt \right\}.\]
When $a >0$, the drift term is bounded from above by $a(r+(a+1)\overline{\sigma}^2/2 ) -
(a+1) \delta $ and we can find $a>0$ such that this bound is negative since $a(r+(a+1)\overline{\sigma}^2/2 ) -
(a+1) \delta \underset{a \rightarrow 0}{\rightarrow} -\delta < 0$. Then, for $x \ge
1$, we have $\Ex[e^{-r \tau^1_x \wedge t }(S^{1}_{\tau^1_x \wedge t})^{1+a}] \le 1$
thanks to Doob's optional sampling theorem. The Fatou lemma gives then $\Ex[e^{-r \tau^1_x
}(S^{1}_{\tau^1_x })^{1+a}] \le 1$, and therefore we get $f_\uparrow(x)=1/\Ex[e^{-r \tau^1_x
}] \ge x^{1+a}$. This shows $\underset{z\rightarrow + \infty}{\lim }   h(z)=0$.
\end{proof}
\begin{proposition}\label{solu0_fg}
  In the case $r=0$, the unique nonincreasing and increasing solution of~(\ref{EDf}) starting
from $1$ in $1$ are respectively:
  \[ f_\downarrow(x)=1, \
  f_\uparrow(x)=\frac{\psi(x)}{\psi(1)}\;\mbox{where}\;\psi(x)= \int_0^x \exp
  \left[\int_1^v   \frac{2 \delta}{u \sigma^2(u)}du \right] dv.\]
  Moreover, we have $f_\uparrow(x)=1/\Px(\tau^1_x <+\infty )$ for $x \ge 1$ and
$f_\uparrow(x)=\Px(\tau^x_1 <+\infty )$ for $x \in (0,1]$.
\end{proposition}
\begin{proof}

When $r=0$, the differential equation $\frac{1}{2} \sigma^2(y) y^2 f''(y)- \delta
y f'(y) =0$ is easy to integrate: $f'(y)=C_3 \exp \left[\int_1^y
  \frac{2 \delta}{u \sigma^2(u)}du \right]$ and then \[g(y)=C_4 + C_3
\int_0^y \exp \left[\int_1^v
  \frac{2\delta }{u\sigma^2(u)}du \right] dv\]
for $C_3,C_4 \in \mathbb{R}.$
For $x \ge 1$,
$f_\uparrow(S^1_{\tau^1_x \wedge t})$ is a bounded martingale that converges almost
surely to $f_\uparrow(x) \mathbf{1}_{\{\tau^1_x <+\infty \}}$ thanks to~(\ref{Dub})
and $\underset{x \rightarrow 0}{\lim }f_\uparrow(x)=0$. Therefore $f_\uparrow(x)\Px(\tau^1_x <+\infty )=1$, and the proof is the same for  $x \in (0,1]$.
\end{proof}


To conclude this section, we state a comparison result which will enables us to
compare $\xss(y)$ and $\ysp_{\sigma}(y)$ with the exercise boundaries obtained in the Black-Scholes model with
constant volatility $\underline{\sigma}$ (resp. $\overline{\sigma}$) where
$\underline{\sigma}$ (resp. $\overline{\sigma}$) bounds the function $\sigma$ from
below (resp. above).
\begin{proposition}\label{prop_compar}
  Let us consider two volatility functions $\sigma_1$ and $\sigma_2$ such that $\forall x>0,\
  \sigma_1(x) \le \sigma_2(x)$  and that satisfy $(\Hs)$. We also assume that either
  $f_{\downarrow,\sigma_1}''$ or $f_{\downarrow,\sigma_2}''$ (resp. either
  $f_{\uparrow,\sigma_1}''$ or $f_{\uparrow,\sigma_2}''$) are nonnegative
  functions   and $r>0$
  (resp. $\delta>0$). Then, we have
  \[  \forall x,y>0, \ P_{\sigma_1}(x,y) \le  P_{\sigma_2}(x,y) \ \ (resp. \
  C_{\sigma_1}(x,y) \le  C_{\sigma_2}(x,y))\]
  and we can compare the exercise boundaries:
  \[ \forall y>0, x^*_{\sigma_1} (y) \ge x^*_{\sigma_2}(y) \ \ (resp. \
  \ysp_{\sigma_1}(y) \le \ysp_{\sigma_2}(y) ) . \]
\end{proposition}
Here and in the proof below, we add in the notation for each mathematical object the volatility function to which
it refers. El Karoui and al. \cite{KJS} and Hobson \cite{Hob} prove that for a convex
payoff function, the price of an American option with finite maturity is a convex function of
the underlying spot price. They deduce monotonicity with respect to the local
volatility function. Their results imply at the same time the convexity assumption
made in the above proposition and its conclusion. In this paper, we prefer to give
autonomous proofs of these results in our simple framework. And we will first use
proposition \ref{prop_compar} to compare with the Black-Scholes case where convexity
is obvious. We can then deduce (Lemma \ref{charac_bound}) that $f''_\downarrow$ and
$f''_\uparrow$ are positive for any $\sigma$ satisfying $(\Hs)$.


\begin{proof}
  Let us consider for example the put case with $f_{\downarrow,\sigma_1}''\ge 0$. Let
  $x \ge z>0$.
  Ito's formula gives:
  \begin{eqnarray*}
    de^{-rt}f_{\downarrow,\sigma_1}(S^{x,\sigma_2}_t)&=&e^{-rt}
    f_{\downarrow,\sigma_1}'(S^{x,\sigma_2}_t) \sigma_2(S^{x,\sigma_2}_t) S^{x,\sigma_2}_t dW_t +e^{-rt}\left[
      \frac{\sigma_2^2(S^{x,\sigma_2}_t)}{2}(S^{x,\sigma_2}_t)^2f''_{\downarrow,\sigma_1}(S^{x,\sigma_2}_t) \right. \\
    && +
      (r-\delta)S^{x,\sigma_2}_t
      f_{\downarrow,\sigma_1}'(S^{x,\sigma_2}_t)-rf_{\downarrow,\sigma_1}(S^{x,\sigma_2}_t) \Big] dt  .
    \end{eqnarray*}
    The term between brackets is nonnegative since we have $ \sigma_1 \le
    \sigma_2$ and $f_{\downarrow,\sigma_1}$ is a convex function solving~(\ref{EDf}). Therefore we get
   $\Ex[ e^{-r \nu_n \wedge
        \tau^{x,\sigma_2}_{z}
      }f_{\downarrow,\sigma_1}(S^{x,\sigma_2}_{\nu_n \wedge
        \tau^{x,\sigma_2}_z } )] \ge
    f_{\downarrow,\sigma_1}(x)$  using the optional sampling theorem, where $\nu_n=\inf \{t \ge 0, S^{x,\sigma_2}_t \ge n
    \} \wedge n$. Since $\nu_n \underset{a.s.}{\rightarrow} +\infty$ and $f_{\downarrow,\sigma_1}(S^{x,\sigma_2}_{\nu_n \wedge
        \tau^{x,\sigma_2}_z } )$ is bounded by
      $f_{\downarrow,\sigma_1}(z)$, Lebesgue's
    dominated convergence theorem gives then 
    $\Ex [ e^{ -r \tau^{x,\sigma_2}_z }]\ge f_{\downarrow,\sigma_1}(x)/
     f_{\downarrow,\sigma_1}(z)$ and so $\Ex [ e^{ -r \tau^{x,\sigma_2}_z }]\ge \Ex [
     e^{ -r \tau^{x,\sigma_1}_z }]$ using Remark~\ref{Tps_atteinte}. The same
     conclusion holds when $f_{\downarrow,\sigma_2}$ is convex by estimating
$\Ex [ e^{ -r \tau^{x,\sigma_1}_{z} } f_{\downarrow,\sigma_2}(S^{x,\sigma_1}_{\tau^{x,\sigma_1}_{z}})]$.
     We then get the result: if  $P_{\sigma_1}(x,y)>(x-y)^+$,
     $P_{\sigma_1}(x,y)=(y-x^*_{\sigma_1}(y))\Ex [ e^{ -r \tau^{x,\sigma_1}_{x^*_{\sigma_1}(y)} }] \le
     (y-x^*_{\sigma_1}(y))\Ex [ e^{ -r \tau^{x,\sigma_2}_{x^*_{\sigma_1}(y)} }]=\Ex [ (y-S^{x,\sigma_2}_{\tau^{x,\sigma_2}_{x^*_{\sigma_1}(y)}})^+ e^{ -r \tau^{x,\sigma_2}_{x^*_{\sigma_1}(y)} }]  \le P_{\sigma_2}(x,y)$.
Now we just observe that $\{x>0, P_{\sigma_1}(x,y) >(x-y)^+
      \} \subset \{x>0, P_{\sigma_2}(x,y) >(x-y)^+ \}$ and thus $x^*_{\sigma_1}
      (y)= \inf \{x>0, P_{\sigma_1}(x,y) >(x-y)^+
      \} \ge \inf \{x>0, P_{\sigma_2}(x,y) >(x-y)^+ \}=x^*_{\sigma_2}(y).$
    \end{proof}

\section{Framework and notations}
We will present in this section the framework that we will consider in all the paper.
To clarify the duality, we will use names that implicitly refer either to the
primal (or ``real'') world, or to the dual world. This denomination has no mathematical
meaning since, as we will see, there are no difference between them. On the contrary,
from a financial point of view, natural variables such as the interest rate, the
dividend rate have their true meaning in the primal world, while in the dual
world they interchange their role. This is the reason why we also name the primal world
``real'' world.\\

\subsubsection*{The primal (``real'') world}

The primal world is the framework we just have described. The spot interest rate $r$
is constant and nonnegative, and $S_t$ is an asset which pays a constant dividend rate $\delta \ge 0$ and is
driven by a homogeneous volatility function $\sigma: \mathbb{R}^*_+ \rightarrow
\mathbb{R}^*_+$ that satisfies  $(\Hs)$ under a risk-neutral measure.  The prices of
the perpetual American put and call are respectively denoted by $P_\sigma (x,y)$ and
$C_\sigma (x,y)$, and their exercise boundary by $\xss$ and $\ysp_{\sigma}$.

\subsubsection*{The dual world}

 In the dual world, $\delta$ plays the role of the interest rate and $r$ of
 the dividend rate; $x$ plays the role of the strike and $y$ is the spot value of the
 share. Let $\eta: \mathbb{R}^*_+ \rightarrow \mathbb{R}^*_+$  be an
 homogeneous volatility function  that is also assumed to satisfy
 $(\Hs)$. We consider then $(\overline{S}^{y}_t,t \ge 0)$ the solution of
 $d\overline{S}_t=\overline{S}_t((\delta-r)dt+\eta(\overline{S}_t)dW_t)$ that starts from $y$ at time $0$. 
 Under that model, we denote respectively by
\[ p_\eta (y,x)= \underset{\tau \in
   \mathcal{T}_{0,\infty}} {\sup} \Ex \left[e^{-\delta \tau}(x-\overline{S}^{y}_{\tau})^+
\right] \text{ and } c_\eta (y,x)= \underset{\tau \in
   \mathcal{T}_{0,\infty}} {\sup} \Ex \left[e^{- \delta \tau}(\overline{S}^{y}_{\tau}-x)^+
 \right]  \]
the prices of the perpetual put and call with strike $x>0$ and spot
 $y$. We can define, as in the primal world,  $g_\downarrow$ and $g_\uparrow$ as the
 unique decreasing (non-increasing when $\delta=0$) and increasing positive solution to
 \begin{equation}\label{EDg}
   \frac{1}{2} \eta^2(x) x^2 g''(x)+ (\delta-r) x  g'(x) -  \delta g(x) =0,
 \end{equation}
 and  we name the exercise boundaries $\xsd_\eta (x)<x$ and $\ysd_\eta (x)>x$ that are
 respectively associated with the put $p_\eta (y,x)$ and the call $ c_\eta (y,x)$.

\subsubsection*{Notations}

 The aim of this paper is to put in evidence a
duality relation and interpret put (resp. call)
prices in the primal world as call (resp. put) prices in the dual
world for a specific volatility function $\eta=\tilde{\sigma}$ (resp.
$\eta=\hat{\sigma}$). When $r=0$ (resp. $\delta=0$), this is trivial since $P_\sigma(x,y)= c_\eta
(y,x)=y$ (resp. $C_\sigma(x,y)= p_\eta (y,x)=x$) but not really fruitful, and we
take thus the following convention in the sequel.
\begin{convention}\label{rdeltapositifs}
  We will always assume $r>0$ (resp. $\delta>0$) to state properties on $P_\sigma$
  and $c_\eta$ (resp.  $C_\sigma$ and $p_\eta$).
\end{convention}

Both worlds being mathematically equivalent, 
{\it we will work with the put price in the primal world and the call price in the dual
world} in order not to do the things twice. Following the Convention
\ref{rdeltapositifs}, we will consider a 
positive interest rate. Let us then denote from now:
\[\boxed{f=f_\downarrow \text{ and } g=g_\uparrow,}\]
and define $\alpha(y)=\frac{y-x^*_\sigma(y)}{f(x^*_\sigma(y))}$ and
$\beta(x)=\frac{\ysd_\eta(x)-x}{g(\ysd_\eta(x))}$. The functions $\alpha$ and  $\beta$ are
positive functions and it follows from the present section that :
\begin{eqnarray}
\forall y>0,\  \forall x \ge x^*_\sigma(y), \   P_\sigma(x,y)
&=& \alpha(y) f(x) \label{Prod_put} \\ 
\forall x>0, \  \forall y \le \ysd_\eta(x), \  c_\eta(y,x) &=&\beta(x)g(y). \label{Prod_call}
\end{eqnarray}
That product form will play an important role for the duality. Let us finally
introduce notations relative the the Black-Scholes model. We define for $\varsigma>0$
 \[a(\varsigma)=\frac{\delta-r+\varsigma^2/2 -
  \sqrt{(\delta-r+\varsigma^2/2)^2 + 2 r \varsigma^2}}{\varsigma^2}<0,\]\[ \  b(\varsigma)=\frac{r-\delta+\varsigma^2/2 +
  \sqrt{(\delta-r-\varsigma^2/2)^2+2\delta \varsigma^2}}{\varsigma^2}=1-a(\varsigma)>1, \]
and we can easily check that $f(x)=x^{a(\varsigma)}$ (resp. $g=x^{b(\varsigma)}$)
when $\sigma(x)=\varsigma$ (resp. $\eta(x)=\varsigma$). In that case, the unique
solution to (\ref{DEF2x*})
(resp. (\ref{DEF2y*})) is:
\begin{equation}\label{boundBS} x^*_{\varsigma}(y)=\frac{a(\varsigma)}{a(\varsigma)-1} y \ \ (resp.\   y^*_{\varsigma}(x)=\frac{b(\varsigma)}{b(\varsigma)-1}x).\end{equation}
With Proposition \ref{prop_compar}, we deduce very useful estimations on the exercise boundaries.

\begin{lemma} 
If, $\forall x>0$,  $\underline{\sigma}\le\sigma(x)\le\overline{\sigma}$
and $\underline{\sigma}\le\eta(x)\le\overline{\sigma}$, then we have~:
\begin{eqnarray}\label{bound_boundaries_x}
\frac{a(\overline{\sigma})}{a(\overline{\sigma})-1}y   \le  \xss(y) \le
\frac{a(\underline{\sigma})}{a(\underline{\sigma})-1}y  & \text{with}& \frac{a(\underline{\sigma})}{a(\underline{\sigma})-1} < \min(1,r/\delta)\\ \label{bound_boundaries_y}
 \frac{b(\underline{\sigma})}{b(\underline{\sigma})-1}x \le
\yse(x) \le \frac{b(\overline{\sigma})}{b(\overline{\sigma})-1}x &\text{with}& \max(1,\delta/r) < \frac{b(\underline{\sigma})}{b(\underline{\sigma})-1}.
\end{eqnarray}
\end{lemma}
\begin{proof}
It is straightforward from Proposition \ref{prop_compar} and (\ref{boundBS}) to get the
bound on the exercise boundaries since $x \mapsto x^\gamma$ is convex for $\gamma
\not \in ]0,1[$. We have to show that for $\varsigma>0$,
$\frac{a(\varsigma)}{a(\varsigma)-1} < \min(1,r/\delta)$ and
$\frac{b(\varsigma)}{b(\varsigma)-1} > \min(1,\delta/r)$. Since $a(\varsigma)<0$ and
$b(\varsigma)>1$, we get that $a(\varsigma)/(a(\varsigma)-1) \in (0,1)$ and
$b(\varsigma)/(b(\varsigma)-1) \in (1,+\infty)$. We can also check that $a(\varsigma)/(a(\varsigma)-1)$ is a root of the polynomial $Q(X)= \delta X^2 -
(r+\delta +\varsigma^2/2) X +r$. As $Q(x)=0 \iff   \sigma^2 x =2(1-x)(r-\delta x)$ and since $a(\varsigma)/(a(\varsigma)-1) \in (0,1)$, we then deduce that
$a(\varsigma)/(a(\varsigma)-1) < r/\delta$.  In the same way, we have $b(\varsigma)/(b(\varsigma)-1) > \delta/r$.\end{proof}

\section{ODE for the exercise boundary}\label{Sec_per_case}

We have seen previously that the exercise boundaries satisfy
\begin{equation}\label{Bound_cond}
  \xss(y)-y=f(\xss(y)) / f'(\xss(y)) \ \  (resp.  \ 
    \yse(x)-x=g(\yse(x))/g'(\yse(x))).
  \end{equation}

We will soon prove
that for fixed $y>0$ (resp. $x>0$), (\ref{Bound_cond}) admits a
  unique solution $\xss(y)$ (resp. $\yse(x)$). 
\begin{lemma}\label{charac_bound}
The function $f'$ (resp. $g'$) is negative (resp. positive) and $f''$ (resp. $g''$)
  is positive on $(0,+\infty)$. Moreover, the boundaries $\xss(y)$ and
  $\yse(x)$ are respectively the unique solution to
  $y-x+f(x)/f'(x)=0$ and $y-x-g(x)/g'(x)=0$. Last, $\xss(y)$,
  $\alpha(y)$, $\yse(x)$ and $\beta(x)$ are $\mathcal{C}^1$ functions on $\mathbb{R}^*_+$.
\end{lemma}
\begin{remark}\label{cvxty}
Positivity of $f''$ and $g''$ and
(\ref{Prod_put}) and~(\ref{Prod_call}) imply positivity of
$\partial^2_xP_\sigma(x,y)$ and $\partial^2_yc_\eta(y,x)$ in the continuation
regions. 
\end{remark}


Differentiating (\ref{Bound_cond}) with respect to $y$ (resp. $x$), one obtains
$1=(\xss)'(y)\frac{f(\xss(y))f''(\xss(y))}{f'(\xss(y))^2}$ (resp.
$1=(\yse)'(y) \frac{g(\yse(y))g''(\yse(y))}{g' (\yse
  (y))^2}$). Using (\ref{Bound_cond}) and equation (\ref{EDf}) (resp. (\ref{EDg}))
one deduces the following result (see equation (\ref{calc_ode}) below).

\begin{proposition}\label{propODE}
Let us assume that the volatility functions $\sigma$ and $\eta$ satisfy $(\Hs)$.
Then, the boundaries $\xss(y)$ and $\yse(x)$ satisfy the following ODEs:
\begin{equation}\label{EDx*}
(\xss)'(y)=\frac{\xss(y)^2\sigma(\xss(y))^2}{2(y-\xss(y))(ry-\delta \xss(y))},
\end{equation} 
\begin{equation}\label{EDy*}(\yse)'(x)= \frac{\eta^2(\yse(x)) \yse(x)^2}{2(\yse(x)-x)(r\yse(x) -
  \delta x)}. 
\end{equation} 
\end{proposition}

\begin{proof}[Proof of Lemma \ref{charac_bound}]
We only give the proof in the put case, the argument being similar for
  the call. By \eqref{EDf}, for $x>0$, $f''(x)$ has the same sign as
  $h(x)=rf(x)+(\delta-r)xf'(x)$. 
If for some $x>0$, $f'(x)=0$, then since $f$ is positive, $f''(x)>0$. Therefore $x$ is a local minimum point
  of
  $f$ which contradicts the decreasing property of this function. Hence
  $f'$ is a negative function.\\
When $\delta\leq r$, $h$ and therefore $f''$ are positive
  functions. When $\delta>r$, we remark that if $f''(x)=0$ then
  $h'(x)=\delta f'(x)<0$. Since the continuous function $f''$ and $h$ have the
  same sign, this implies that 
\begin{equation}
   \forall x >
  \inf\{z>0:f''(z)\leq 0\},\;f''(x) < 0.\label{borninff}
\end{equation}Now for $y>0$, by~(\ref{EDf}) then~(\ref{Bound_cond}) , we have
  \begin{equation}\label{calc_ode}\frac{\xss(y)^2\sigma(\xss(y))^2}{2}\frac{f''(\xss(y))}{f'(\xss(y))}=r\frac{f(\xss(y))}{f'(\xss(y))}-(r-\delta)\xss(y) = \delta \xss(y)-ry. \end{equation}
  By \eqref{bound_boundaries_x}, the right-hand-side is negative and moreover
  $\lim_{y\rightarrow +\infty}\xss(y)=+\infty$. Hence
$\sup\{z>0:f''(z)>0\}=+\infty$ and with \eqref{borninff}, we conclude that $f''>0$.


According to~(\ref{Bound_cond}), $F(\xss(y),y)=0$ where \[F(x,y)=y-x+f(x)/f'(x).\]
The function $F$ is $\mathcal{C}^1$ on $(0,+\infty)\times(0,+\infty)$ and such that 
\[\forall x,y>0,\;\partial_xF(x,y)=-f(x)f''(x)/f'
  (x)^2<0.\] Therefore for fixed $y>0$, $x^*(y)$ is the unique solution
  to $F(x,y)=0$. Moreover, $y\rightarrow x^*(y)$ is $C^1$ by the
  implicit function theorem. Last, one deduces from \eqref{Bound_cond}
  that $\alpha(y)$ is a $C^1$ function.
\end{proof}

The positivity of $f''$ and $g''$ gives the following result.
\begin{corollary}\label{cor_compar} The comparison result stated in Proposition \ref{prop_compar} holds for any
  $\sigma_1 \le \sigma_2$ satisfying $(\Hs)$.
\end{corollary}


Let us now give a uniqueness result for the ODEs~(\ref{EDx*}) and~(\ref{EDy*}).

\begin{proposition}\label{Unicity}
  There is only one solution $\xss$ of~(\ref{EDx*}) (resp. $\yse$ of~(\ref{EDy*}))
  defined on $(0,+\infty)$ that satisfies
  $\forall y>0,\;c_1 y \le \xss(y) \le c_2 y$ with $0<c_1 \le c_2 < \min(1,r/\delta)$ (resp. $\forall x>0,\;d_1 x \le
  \yse(x) \le d_2 x$ with $d_1> \max(1,\delta/r)$).
\end{proposition}
\begin{proof}
  Let us first remark that the uniqueness result for~(\ref{EDx*}) is equivalent to the uniqueness
  result for~(\ref{EDy*}). Indeed, it is easy to see that $\xss(y)$ is solution of~(\ref{EDx*}) if and only if $\hat{y}(x):=1/(\xss(1/x))$ is solution of~(\ref{EDy*})
  with the volatility function $\eta(x)=\sigma(1/x)$. This new volatility also
  satisfies $(\Hs)$. Moreover, $d_1 x \le \hat{y}(x) \le d_2 x$ with $d_1 > \max(1,\delta/r)$ if, and only if 
  $0 \le c_1 y \le \xss(y) \le c_2 y$ with $0<c_1 \le c_2 < \min(1,r/\delta)$.
  
  Let us suppose then that there are two solutions of~(\ref{EDy*}),  $y_1(x)$ and
  $y_2(x)$, that are defined on $\mathbb{R}_+$ and satisfy   $d_2 x \ge y_{j}(x) \ge d_1 x$
  for some $d_2>d_1> \max(1,\delta/r)$. Since $y_{1}'(x)>0$ for $x>0$, $y_1$ is invertible
  and we have:
  \begin{eqnarray*}
    \frac{d}{dx}y_1^{-1}(y_2(x)) &=& \frac{ y_2(x)^2 \eta(y_2(x))^2 }{2(y_2(x)-x)(r y_2(x)-
      \delta x) }\frac{2(y_2(x)-y_1^{-1}(y_2(x)))(ry_2(x)- \delta y_1^{-1}(y_2(x))
      )}{y_2(x)^2 \eta(y_2(x))^2 } \\
    &=& \frac{(y_2(x)-y_1^{-1}(y_2(x)))(ry_2(x)- \delta y_1^{-1}(y_2(x))}{(y_2(x)-x)(r y_2(x)-
      \delta x)}.
  \end{eqnarray*}
Thus, the function $\psi(x)=y_1^{-1}(y_2(x))/x$  solves
\begin{eqnarray}\label{eqpsi}
\psi'(x)&=&\frac{1}{x}\left[\frac{y_2(x)-\psi(x)x}{y_2(x)-x} \times \frac{ry_2(x) - \delta x
    \psi(x) }{ry_2(x)- \delta  x} -\psi(x) \right] \nonumber\\
&=& \frac{1}{x}\left[\left(1- \frac{\psi(x)-1}{y_2(x)/x-1} \right) \left( 1-
    \frac{
    \psi(x) - 1}{ry_2(x)/( \delta  x)-1}\right) -\psi(x) \right] . 
\end{eqnarray}
The estimation  $d_2 x \ge y_{j}(x) \ge d_1 x$ for $j \in \{ 1,2\}$ with
$d_1>\max(1,\delta/r)$ implies that :
\begin{align}\label{CONDpsi1}
  &\exists A>0,\forall x>0, 1/A \le \psi(x) \le A,\\
&\forall
x>0,\;\psi(x)<\min\left(\frac{y_2(x)}{x},\frac{ry_2(x)}{\delta x}\right),\;\frac{y_2(x)}{x}-1>0\text{ and }  \frac{ry_2(x)}{\delta  x}-1>0\label{CONDpsi2}.  \end{align}

Since local uniqueness holds for
\eqref{eqpsi} by the Cauchy Lipschitz theorem, the only solution
$\varphi$ such
that $\varphi(1)=1$ is the constant $\varphi \equiv 1$. Therefore
checking that \eqref{CONDpsi1} does not hold for solutions $\varphi$ satisfying
\eqref{CONDpsi2} and such that
$\varphi(1)\neq 1$ is enough to conclude that $\psi\equiv 1$.

Let $\varphi$ be a solution to \eqref{eqpsi} satisfying
\eqref{CONDpsi2}. If $\varphi(1)>1$, by local uniqueness for
\eqref{eqpsi}, for all $x \in \mathbb{R}_+^*$, $\varphi(x)>1$. By
\eqref{CONDpsi2}, one deduces that for all $x \in \mathbb{R}_+^*$,
$\varphi'(x)<\frac{1-\varphi(x)}{x}<0$. Therefore, $\varphi' (x) \le
(1-\varphi(1))/x$ for $x \in (0,1]$, and we have
\[ \varphi(x) \ge \varphi(1)+(1-\varphi(1)) \ln(x) \underset{x \rightarrow 0}{\rightarrow}
+\infty \] which is contradictory to \eqref{CONDpsi1}. In the same manner, if $\varphi(1)<1$, $\varphi(x)<1$ for
$x \in \mathbb{R}_+^*$ and $\varphi$ is strictly increasing. In particular, for $x \le 1$,
$\varphi'(x) \ge (1-\varphi(1))/x$ and therefore $\varphi(1)-\varphi(x) \ge (1-\varphi(1)) \ln (1/x)
\underset{x \rightarrow 0}{\rightarrow} +\infty$ and this yields another
contradiction.\end{proof}

\begin{corollary}\label{sig->x*}
Let us denote $\tilde{\mathcal{C}}=\{f \in
  \mathcal{C}^1(\mathbb{R}_+^*),  \text{ s.t. } f(0)=0,\  \exists 0<a<b,\forall  x \ge 0,
  \  a\le f'(x) \le b\}$.
  The application $\sigma \mapsto \xss$ (resp. $\eta \mapsto \yse$) is 
one-to-one between the set $\{ \sigma \in \mathcal{C}(\mathbb{R}_+^*) \text{ that
  satisfies } (\Hs) \}$ and the set of function $\tilde{\mathcal{C}}_x=\{x \in \tilde{\mathcal{C}}
  ,  \text{ s.t. } \exists 0<c_1 \le c_2  <
  \min(1,r/\delta),\forall y>0, \  c_1y \le x(y) \le c_2 y \}$ (resp. $\tilde{\mathcal{C}}_y=\{y \in \tilde{\mathcal{C}}
 ,  \text{ s.t. }\exists \max(1,\delta/r)<d_1 \le d_2  ,\forall x>0, \  d_1 x \le y(x) \le d_2 x\}$.)
\end{corollary}
\begin{proof}
  If $\sigma$ is a continuous function satisfying $(\Hs)$, by \eqref{EDx*} and \eqref{bound_boundaries_x}, $\xss$ belongs to $\tilde{\mathcal{C}}_x$.
The one to one property is easy to get. If $x^*_{\sigma_1}\equiv x^*_{\sigma_2}$ with 
  $\sigma_1$ and $\sigma_2$ satisfying $(\Hs)$, the ODE~(\ref{EDx*}) ensures that
  $\sigma_1^2(x^*_{\sigma_1}(y))=\sigma_2^2(x^*_{\sigma_2}(y))$ for $y>0$. Therefore
  $\sigma_1\equiv\sigma_2$.\\
  Let us check the
  onto property and consider $x^*(y) \in \tilde{\mathcal{C}}_x$. The function $\sigma$ defined by
  \begin{equation}\label{sigma_de_x*}
    \sigma(x^*(y))=\frac{\sqrt{2(y-x^*(y))(ry-\delta x^*(y)){x^*}'(y) }}{x^*(y)} \end{equation}
  is well defined thanks to the hypothesis made on $x^*$. As $\xss$
  satisfies \eqref{bound_boundaries_x} and solves the
  same ODE~(\ref{EDx*}) as $x^*$, we have $x^*\equiv\xss$ using Proposition
  \ref{Unicity}.\\The proof for $\eta \mapsto \yse$ is the same and gives incidentally
  the expression of $\eta$ in function of the exercise boundary $y^*(x)$:
  \begin{equation}\label{eta_de_y*}
    \eta(y^*(x))=\frac{\sqrt{2(y^*(x)-x)(ry^*(x)-\delta x){y^*}'(x) }}{y^*(x)}. \end{equation}
  \end{proof}

\section{The call-put duality}\label{Sec_dual}

This section is devoted to the key result of the paper : for related
local volatility functions $\sigma$ and $\eta$, we can interpret
a put price in the primal world as a call price in the dual world.

\subsection{The main result}
\begin{theorem}[Duality]\label{DualProp}
  The following conditions are equivalent:
  \begin{enumerate}
    \item
\begin{equation}\label{Dual} \forall x,y>0,\  P_\sigma(x,y)= c_\eta(y,x). 
\end{equation}
\item $\xss$ and $\yse$ are reciprocal functions: $\forall x>0, \ \xss(\yse(x))=x$.
\item $\eta\equiv\tsigma$ where
\begin{equation}\label{Dual_sig} \tsigma(y)= \frac{2(y-\xss(y))(ry-\delta \xss(y))}{y \xss(y)
    \sigma (\xss(y))}.\end{equation}
\item $\sigma\equiv\uteta$ where \begin{equation}\label{Dual_eta} \uteta(x)=
    \frac{2(\yse(x)-x)(r \yse(x)-\delta x)}{ \yse(x) x
    \eta (\yse(x))}.\end{equation}
  \end{enumerate}
\end{theorem}

\begin{remark}
  Thanks to relation~(\ref{bound_boundaries_x}) (resp.~(\ref{bound_boundaries_y})), if $\sigma$ (resp. $\eta$) satisfies $(\Hs)$ then the dual volatility function $\tsigma$ defined
  by~(\ref{Dual_sig}) (resp. $\uteta$ defined
  by~(\ref{Dual_eta})) satisfies $(\Hs)$.
\end{remark}

\begin{proof}
$1 \implies 2$~: We have on the one hand $P_\sigma(x,y)=y-x$ on $\{(x,y), \ x \le \xss(y)
\}$ and $P_\sigma(x,y)>y-x$ outside, and on the other
hand $c_\eta(y,x)=y-x$ on $\{(x,y), \ y \ge \yse(x)
\}$ and $c_\eta(y,x)>y-x$ outside. The duality
relation~(\ref{Dual}) imposes then that $\{(x,y), \ x \le \xss(y)
\}=\{(x,y), \ y \ge \yse(x)
\}$ and so $\yse (\xss (y))=y$.

$2 \implies 3,4$~: Taking the derivative of the last relation, we get thanks to
(\ref{EDx*}) and~(\ref{EDy*})\\
$\frac{\xss (y)^2\sigma(\xss (y))^2}{2(y-\xss (y))(ry-\delta \xss(y))} \frac{\eta^2(y) y^2}{2(y-\xss(y))(ry -
  \delta \xss(y))}=1$ and deduce~(\ref{Dual_sig}) and~(\ref{Dual_eta}).

$3 \implies 2$~ (resp. $4 \implies 2$~) : By~(\ref{EDx*})
(resp.~(\ref{EDy*})) and~(\ref{Dual_sig}) (resp.~(\ref{Dual_eta})), ${\xss}^{-1}$ (resp. ${\yse}^{-1}$)
satisfies~(\ref{EDy*}) (resp.~(\ref{EDx*})). Since by
(\ref{bound_boundaries_x}) (resp.~(\ref{bound_boundaries_y})) this
function satisfies~(\ref{bound_boundaries_y})
(resp.~(\ref{bound_boundaries_x})), one concludes by
Proposition \ref{Unicity}.

$2 \implies 1$~: The equality (\ref{Dual}) is clear in the exercise region since
$\{(x,y),x \le \xss(y) \} =\{(x,y),y \ge \yse(x) \}$. Let us check that it also holds in
the continuation region. Using the product form~(\ref{Prod_call}), and the smooth-fit principle (Theorem \ref{Opt_stop_call}) we get for all $y\in \mathbb{R}^*_+$
\begin{equation*}
 \left\{  \begin{array}{l}
y-\xss(y)=\beta(\xss(y))g(y) \\
1=-\beta(\xss(y))g'(y).
\end{array}
 \right.   \end{equation*}
Differentiating the first equality with respect to $y$, one gets $1-\xss(y)'=\xss(y)'\beta'(\xss(y))g(y) +
\beta(\xss(y))g'(y)$, which combined with the second equality gives
\[-1=\beta'(\xss(y))g(y).\]
Dividing by the first equality and using (\ref{Bound_cond}), one deduces
$\frac{\beta'}{\beta}(\xss(y))=\frac{f'}{f}(\xss(y))$. Since $\xss: \mathbb{R}_+^*
\rightarrow  \mathbb{R}_+^*$ is a bijection, there is a constant $C \not =0$ such
that $\beta \equiv Cf$.
Since $\forall y>0, \alpha(y)f(\xss(y))=y-\xss(y)=\beta(\xss(y))g(y)$, one has
$\alpha  \equiv g/C$. From~(\ref{Prod_put}) and~(\ref{Prod_call}), one concludes
that~(\ref{Dual}) holds. \end{proof}

In this proof, we have shown that $\alpha$ is proportional to $g$, and so there is a
constant $C>0$ such that $P_\sigma(x,y)=Cf(x)g(y)$ for $x \le \xss(y)$.
In the Black-Scholes'~case, we have $P_\sigma(x,y)=Cx^{a(\sigma)}y^{b(\sigma)}$ for $x \le \xss(y)$, and we are able to calculate
$C$ using the boundary condition $P_\sigma(\xss(y),y)=y-\xss(y)$. We
get that $C=b(\sigma)^{-b(\sigma)}/(-a(\sigma))^{a(\sigma)}$. We retrieve then the
already known analytical formulae for the put prices (e.g. Gerber and Shiu \cite{GS}).


\subsection{An analytic example of dual volatility functions }\label{an_ex}
By \eqref{Dual_eta} and \eqref{EDy*}, if $y^*\in\tilde{\mathcal C}_y$
(where $\tilde{\mathcal C}_y$ is defined in Corollary \ref{sig->x*}),
then the reciprocal function of $y^*$ is the put exercise
boundary $x^*_\sigma$ associated to the local volatility function 
$$\sigma(x)=\frac{\sqrt{2(ry^*(x)- \delta x)(y^*(x)-x)}}{x
  \sqrt{y^*(x)'}}.$$
Now by \eqref{eta_de_y*}, $y^*$ is the call
  exercise boundary associated with the dual volatility function :

$$\tilde{\sigma}(y)=\frac{\sqrt{2(y-\xss(y))(ry-\delta \xss(y)){y^*}'(\xss(y)) }}{y}.$$

Let us consider the family of exercise boundaries \[
y^*(x)=x\frac{x+a}{bx+c} \]
where $a,b,c$ are positive constants such that
$\max(c/a,b)<\min(1,r/\delta)$ (condition ensuring $y^* \in \tilde{\mathcal C}_y$).
Since $y^*(x)'=(bx^2+2cx+ac)/(bx+c)^2$, 
one has \[\sigma(x)=\sqrt{ 2 \frac{((r-\delta b)x +ra - \delta c )((1- b)x + a - c
    )}{bx^2+2cx+ac}}, x> 0.\]
Moreover, the function $\xss(y)$ is the only positive root of
the polynomial function: $X^2+X(a-by)-cy$, that is:
\[\xss(y)=\frac{1}{2}\left( by-a+\sqrt{(by-a)^2+4cy } \right)\]
and \[ \forall y >0, \ 
\tilde{\sigma}(y)=\frac{\sqrt{2(y-\xss(y))(ry-\delta \xss(y) )(b\xss(y)^2+2c\xss(y)+ac) }}{y(b \xss(y) +c)}. \]
\begin{figure}[h]
  \psfrag{psigtilde}{ $p_{\tsigma}$}
  \psfrag{csigtilde}{ $c_{\tsigma}$ }
  \psfrag{Psigma}{ $P_\sigma$ } 
  \psfrag{Csigma}{ $C_\sigma$ }

  \centerline{ \psfig{file=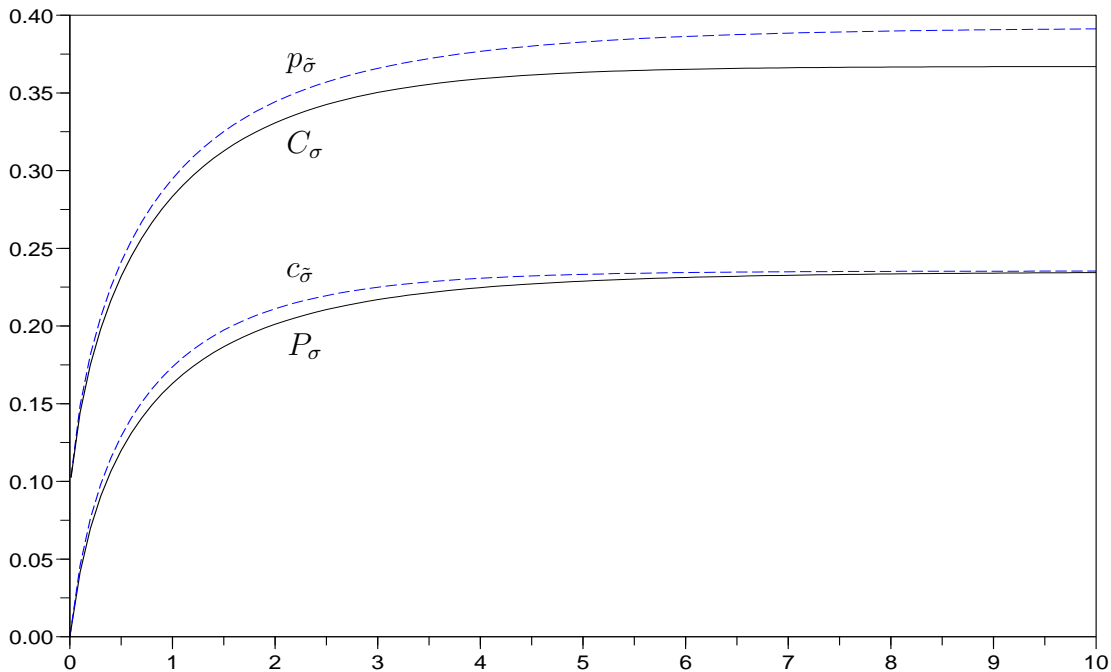,width=18cm,height=11cm} }
    \caption{{\small $P_\sigma(T,x,y)$ and $c_{\tsigma}(T,y,x)$, and $C_\sigma(T,x,y)$
      and $p_{\tsigma}(T,y,x)$ as functions of the maturity $T$ for $x=0.5$, 
      $y=0.4$, $r=0.2$, $\delta=0.1$ and the volatility parameters $(a,b,c)=(1,0.4,0.1)$.}} \label{Dual_inf}
\end{figure} 
This example enables us to check numerically the duality. We have
plotted in Figure~\ref{Dual_inf}, the
prices of an American put $P_\sigma(T,x,y)$ in the primal world for the local volatility $\sigma(x)$ and
an American call $c_{\tsigma} (T,y,x)$ in the dual world for the local volatility 
$\tsigma(x)$ as functions of the maturity $T$. These prices have been computed using a finite difference
method. We can see at $T=10$ that the limit value is quite 
reached and both prices are equal. The plots are nonetheless distinct which means
that the duality does not hold for finite maturities. We have also plotted, in function of $T$, $C_\sigma(T,x,y)$
in the primal world and $p_{\tsigma}(T,y,x)$ in the dual world to check numerically whether the volatility function
$\tsigma$ is such as $C_\sigma(x,y)=p_{\tsigma}(y,x)$. As we can see, the curves do not
seem to converge toward the same limit when $T$ is large. This means
that the volatility function $\hat{\sigma}$ such that $\forall
x,y>0,\;C_\sigma(x,y)=p_{\hat{\sigma}}(y,x)$ (obtained from $\sigma$ as
$\uteta$ is obtained from $\eta$ but with exchange of $r$ and $\delta$) is different from $\tilde{\sigma}$.

\section{ Consequence of the duality~: A (theoretical) method of calibration
  for the volatility  $\sigma(x)$}

In that section, we will put in evidence the
importance of the duality within the calibration scope. 
We suppose for this that we are on a (virtual) market where are traded perpetual
securities, and where the short interest rate $r$ and the dividend rate
$\delta$ can be observed. This means that we know the price of the share
$x_0$, and the market quotes on that share the
perpetual American puts and calls for all strikes
$K>0$. We name respectively $p(K)$ and $c(K)$ these prices and denote:
\begin{equation}\label{XY}  
X=\sup \{K>0, \ c(K)=x_0-K \} \text{ and } Y=\inf \{K>0, \ p(K)=K-x_0 \}.  
\end{equation}
We will first suppose that the put and call prices derive from a
time-homogeneous local volatility model before relaxing this assumption.
\subsection{The calibration procedure}
Let us assume that there is a
volatility function $\sigma$ satisfying ($\Hs$) such that for all $K>0$, 
$p(K)=P_{\sigma}(x_0,K)$ and $c(K)=C_{\sigma}(x_0,K)$. The following
proposition says that these prices characterize $\sigma$ and
its proof gives a constructive way to retrieve the volatility function from the prices.
\begin{proposition}\label{calib}
Let us consider $x_0>0$. The map \[\sigma  \mapsto (
(P_{\sigma}(x_0,K),C_{\sigma}(x_0,K)),K>0)\]  is one-to-one on the set of volatility
functions satisfying $(\Hs)$.
\end{proposition}
\begin{proof}
We first consider the put case. The differential equation satisfied
by the put prices in the continuation region makes only appear 
the values and the derivatives in $x$, $K$ being fixed. Hence, we cannot
exploit directly the prices. But the
duality relation enables to get a differential equation in the strike variable. Thanks to the
Duality Theorem, we have $P_{\sigma}(x_0,K)=c_{\tsigma}(K,x_0)$ for some $\tsigma$ satisfying ($\Hs$). It is then easy to
calibrate $\tsigma(.)$. Indeed, one has $ \frac{K^2
  \tsigma(K)^2}{2}p''(K) + K (\delta-r) 
p'(K) - \delta p(K) =0$ for $K <
Y=\ysd_{\tsigma}(x_0)$. Since the differential equation is valid only for 
$K < Y$, we only get $\tsigma$ on $(0,Y]$ by continuity: \[ \forall
K \le Y, \
\tsigma(K)=\frac{1}{K} \sqrt{\frac{2(\delta p(K) +  K (r-\delta)
    p'(K))}{p''(K)}}\]
which is well defined since $p''(K)=\partial_K^2 c_{\tsigma}(K,x)>0$
(Remark \ref{cvxty}). Then, we can calculate
the exercise boundary $\ysd_{\tsigma}(x)$, for $x\in(0,x_0]$,
solving~(\ref{EDy*}) supplemented with the final condition $ \ysd_{\tsigma}(x_0)=Y$ backward.
This step only requires the knowledge of $\tsigma$ only on the interval $(0,Y]$.
Finally, we can recover the desired volatility $\sigma(x)$ for $x \le
x_0$ thanks to~(\ref{Dual_sig}):
\begin{equation}\label{calib_put}
\forall x \in (0,x_0 ], \
\sigma(x)=\frac{2(\ysd_{\tsigma}(x)-x)(r\ysd_{\tsigma}(x)-
  \delta x) }{x\ysd_{\tsigma}(x)
  \tsigma(\ysd_{\tsigma}(x))} .
\end{equation}

Now let us consider the calibration to the call prices. This relies on the same
principle, but we have to be careful because the Duality Theorem is stated given to
the call interest rate $\delta$ and dividend rate $r$. So we have to interchange
these variables when we apply that theorem.
There is a function $\hat{\sigma}$ satisfying
$(\Hs)$ such that~: $\forall K>0, \
C_{\sigma}(x_0,K)=p_{\hat{\sigma}} (K,x_0)$. We have
\[\frac{1}{2}K^2 \hat{\sigma} (x)^2 c''(K)
+ (\delta-r)K c'(K) - \delta c(K) =0 \] 
for $K>X=\xsd_{\hat{\sigma}}(x_0)$. 
 Thus, we get
\[ \forall K \ge X, \  \hat{\sigma}(K)=\frac{1}{K}
\sqrt{\frac{2 (\delta c(K) +  K (r-\delta)
    c'(K))}{c''(K)}}\]
which is well defined for analogous reasons. We can then obtain as before the exercise
boundary solving~(\ref{EDx*}) forward 
\[\forall y \ge x_0, \
\xsd_{\hat{\sigma}}(y)'=\frac{\xsd_{\hat{\sigma}}(y)^2 \hat{\sigma}(\xsd_{\hat{\sigma}}(y))^2}{2(y-\xsd_{\hat{\sigma}}(y))(\delta y-r
  \xsd_{\hat{\sigma}}(y))}  ,\ \xsd_{\hat{\sigma}}(x_0)=X\]
and we finally get the volatility $\sigma(y)$ for $y \ge x_0$ using the Duality
Theorem. More precisely, we interchange $r$ and $\delta$ in \eqref{Dual_sig} to get
\begin{equation}\label{calib_call}\sigma(y)=\frac{2(y-\xsd_{\hat{\sigma}}(y))(\delta
  y - r
  \xsd_{\hat{\sigma}}(y))}{y
  \xsd_{\hat{\sigma}}(y)
  \hat{\sigma}(\xsd_{\hat{\sigma}}(y)) } .\end{equation}\end{proof}


 


This calibration method, although being
theoretical, sheds light on a striking and interesting result: the
perpetual American put prices only
give the restriction of $\sigma(x)$ to $(0,x_0]$ and the call prices
only the restriction of $\sigma(x)$ to $[x_0,+\infty)$. This has the
following economical interpretation : long-term American put prices mainly give information on the downward
volatility while long-term American call prices give information on the
upward volatility. This dichotomy is remarkable. In comparison,
according to Dupire's
formula \cite{Dup}, there
is no such phenomenon for European options : the knowledge of the call prices gives
the whole local volatility surface, not only one part. In other words, the European
call and put prices give the same information on the volatility while the perpetual American call and
put prices give complementary information.

Thus, one may think that the perpetual American call and put
prices only depend on a part of the volatility curve. This is precised by the Proposition
below that gives necessary and sufficient conditions on the volatility functions to
observe the same put prices (resp. call prices).
\begin{proposition}\label{equivalences}
Let us consider $x_0>0$ and $\sigma_1(.)$, $\sigma_2(.)$ two volatility functions
satisfying ($\Hs$). Then, the following properties are equivalent:
\begin{enumerate}[(i)]
\item $\forall y>0,\  P_{\sigma_1}(x_0,y)=  P_{\sigma_2}(x_0,y)$ (resp. $\forall
  y>0,\  C_{\sigma_1}(x_0,y)= C_{\sigma_2}(x_0,y)$)

\item $\forall y \le \ysd_{\tsigma_2}(x_0), \tsigma_1(y)=\tsigma_2(y)$. (resp.
  $\forall x \ge \xsd_{\hat{\sigma_1}}(x_0),
  \hat{\sigma_1}(x)=\hat{\sigma_2}(x)$ where
 $\hat{\sigma_j}$ denotes the local volatility function such that
 $\forall x,y>0, \;C_{\sigma_j}(x,y)=p_{\hat{\sigma_j}}(y,x)$.)
\item $\forall x\in (0, x_0], \  \sigma_1(x)=\sigma_2(x)$ and $\ysd_{\tsigma_1}(x_0)=\ysd_{\tsigma_2}(x_0) $. (resp.
 $\forall x \in [x_0,+\infty), \  \sigma_1(x)=\sigma_2(x)$ and
 $\xsd_{\hat{\sigma_1}}(x_0)=\xsd_{\hat{\sigma_2}}(x_0)$.)
\item $\forall x \in (0, x_0], \  \sigma_1(x)=\sigma_2(x)$ and
  $\frac{f'_{\downarrow,\sigma_1}(x_0)}{f_{\downarrow,\sigma_1}(x_0)}=\frac{f'_{\downarrow,\sigma_2}(x_0)}{f_{\downarrow,\sigma_2}(x_0)}$. (resp.  $\forall x \in [x_0,+\infty), \  \sigma_1(x)=\sigma_2(x)$ and $\frac{f'_{\uparrow,\sigma_1}(x_0)}{f_{\uparrow,\sigma_1}(x_0)}=\frac{f'_{\uparrow,\sigma_2}(x_0)}{f_{\uparrow,\sigma_2}(x_0)}$.)
  \item $f_{\downarrow,\sigma_1}$ and $f_{\downarrow,\sigma_2}$ (resp. $f_{\uparrow,\sigma_1}$ and $f_{\uparrow,\sigma_2}$) are
    proportional on $(0,x_0]$ (resp. $[x_0,+\infty)$).
  \item  $\forall  x\le x_0, \forall y>0,\  P_{\sigma_1}(x,y)=  P_{\sigma_2}(x,y)$ (resp. $\forall x \ge x_0,\forall
    y>0,\  C_{\sigma_1}(x,y)= C_{\sigma_2}(x,y)$).
    \end{enumerate}
\end{proposition}
\begin{remark}
   \begin{itemize}
\item Among these many conditions, let us remark that condition $(ii)$ on the dual
volatility is much simpler than condition $(iii)$ on the primal
volatility since the latter
requires the equality of the dual exercise boundaries at $x_0$. 
      \item When $\delta=0$, according to Remark \ref{solu_delta0}, in
      the put case, 
      condition $(iv)$ also writes $\forall x\in (0, x_0], \
      \sigma_1(x)=\sigma_2(x)$ and 
   \[{\int^{+\infty}_{x_0}\left( \frac{1}{v^2} \exp \left[- \int_{x_0}^v
      \frac{2r}{u\sigma_1^2(u)}du\right] \right) dv}={ \int^{+\infty}_{x_0}\left( \frac{1}{v^2} \exp \left[- \int_{x_0}^v
        \frac{2r}{u\sigma_2^2(u)}du\right] \right) dv} . \]
\item Since, by definition of $f_{\downarrow,\sigma_j}$
  (resp. $f_{\uparrow,\sigma_j}$) and the strong Markov property
  $\forall 0<z\leq
  x,\;\E[e^{-r\tau^x_{\sigma_j,z}}]=f_{\downarrow,\sigma_j}(x)/
  f_{\downarrow,\sigma_j}(z)$ (resp. $\forall 0<x\leq
  z,\;\E[e^{-r\tau^x_{\sigma_j,z}}]=f_{\uparrow,\sigma_j}(x)/
  f_{\uparrow,\sigma_j}(z)$), the probabilistic counterpart of
  assertion $(v)$ is $\forall 0<z\leq
  x\leq
  x_0,\;\E[e^{-r\tau^x_{\sigma_1,z}}]=\E[e^{-r\tau^x_{\sigma_2,z}}]$
  (resp. $\forall x_0\leq
  x\leq
  z,\;\E[e^{-r\tau^x_{\sigma_1,z}}]=\E[e^{-r\tau^x_{\sigma_2,z}}]$).\end{itemize}
\end{remark}

\begin{proof} We consider for example the put case.

$(i) \implies (ii)$~: See the proof of Theorem \ref{calib}.

$(ii) \implies (iii)$~: Let us define
$\psi(x)=(\ysd_{\tsigma_1})^{-1}(\ysd_{\tsigma_2}(x))/x$. We can show as in the proof
of Proposition \ref{Unicity} that $\psi(x_0)=1$ and then $\psi \equiv 1$ on $(0,x_0]$, otherwise
it would go to $0$ or $+\infty$ when $x \rightarrow 0$, which is not possible thanks
to~(\ref{bound_boundaries_y}). We get then $\forall x\in (0, x_0], \
\sigma_1(x)=\sigma_2(x)$ using~(\ref{Dual_eta}) that express $\sigma_j$ in function of
$\ysd_{\tsigma_j}$ and $\tsigma_j$, $j \in \{1,2\}$.

$(iii) \implies (iv)$~: Thanks to~(\ref{Bound_cond}) and Theorem \ref{DualProp}, we have
$\frac{f'_{\downarrow,\sigma_1}(x_0)}{f_{\downarrow,\sigma_1}(x_0)}=\frac{-1}{\ysd_{\tsigma_1}(x_0)-x_0}=\frac{-1}{\ysd_{\tsigma_2}(x_0)-x_0} =\frac{f'_{\downarrow,\sigma_2}(x_0)}{f_{\downarrow,\sigma_2}(x_0)}$.

$(iv) \implies (v)$~: 
The set of solutions to $\frac{1}{2} \sigma_1^2(x) x^2 f''(x)+ (r-\delta) x  f'(x) -
 rf(x) =0$ on $(0,x_0]$ is a two-dimensional vector space, but thanks to the
 relation $\frac{f'_{\downarrow,\sigma_1}(x_0)}{f_{\downarrow,\sigma_1}(x_0)}=\frac{f'_{\downarrow,\sigma_2}(x_0)}{f_{\downarrow,\sigma_2}(x_0)}$,
 $f_{\downarrow,\sigma_1}$ and $f_{\downarrow,\sigma_2}$ are proportional on $(0,x_0]$~:
\begin{equation}\label{proportion_f} \forall x \le
x_0, \
f_{\downarrow,\sigma_1}(x)=\frac{f_{\downarrow,\sigma_1}(x_0)}{f_{\downarrow,\sigma_2}(x_0)}f_{\downarrow,\sigma_2}(x).
\end{equation}
 
$(v) \implies (vi)$~:
The proportionality implies that $\forall x \in (0, x_0], \
\frac{f_{\downarrow,\sigma_1}(x)'}{f_{\downarrow,\sigma_1}(x)}=\frac{f_{\downarrow,\sigma_2}(x)'}{f_{\downarrow,\sigma_2}(x)}$,
and then  $(\ysd_{\tsigma_1}(x)-x )^{-1}=(\ysd_{\tsigma_2}(x)-x)^{-1}$ using
(\ref{Bound_cond}) and Theorem \ref{DualProp}. Therefore
\[\forall x \in (0, x_0], \ 
\ysd_{\tsigma_1}(x)=\ysd_{\tsigma_2}(x).\]
 We have
$\alpha_{\sigma_1}(\ysd_{{\tsigma}_1}(x) ) 
f_{\downarrow,\sigma_1}(x)= \alpha_{\sigma_2}(\ysd_{\tsigma_2}(x))
f_{\downarrow,\sigma_2}(x)$ using~(\ref{Prod_put}), and obtain from~(\ref{proportion_f}) that
\begin{equation}\label{proportion_alpha} \forall x\leq x_0,\;\forall y \le
\ysd_{\tsigma_1}(x_0), \
\alpha_{\sigma_1}(y)=\frac{f_{\downarrow,\sigma_2}(x_0)}{f_{\downarrow,\sigma_1}(x_0)}
\alpha_{\sigma_2}(y)=\frac{f_{\downarrow,\sigma_2}(x)}{f_{\downarrow,\sigma_1}(x)}
\alpha_{\sigma_2}(y). \end{equation}
Thus, we deduce from~(\ref{Prod_put}),~(\ref{proportion_f}) and~(\ref{proportion_alpha}) the equality of the put prices for the low strikes
\[\forall x\leq x_0,\;\forall y \le \ysd_{\tsigma_1}(x), P_{\sigma_1}(x,y)= P_{\sigma_2}(x,y). \]
For  $y>\ysd_{\tsigma_1}(x)=\ysd_{\tsigma_2}(x)$, the
equality is clear since both prices are equal to $y-x$.

$(vi) \implies (i)$~: clear.\end{proof}

Let us observe that the point $(ii)$ of the last proposition allows to
exhibit different volatility functions with analytic expressions that give the same put (or call)
prices. Let us consider the same family as in subsection \ref{an_ex} coming from
the call exercise boundary $y^*_1(x)=x\frac{x+a}{bx+c}$ (assuming $a,b,c>0$ and $\max(c/a,b)<\min(1,r/\delta)$). For
$x_0>0$, we introduce the exercise boundary:
\[y^*_2(x)=y^*_1(x) \text{ for } x \le x_0 \text{ and }
y^*_2(x)=y^*_1(x_0)+(y^*_1)'(x_0)(x-x_0) \text{ for } x \ge x_0.\]
 The condition ($y^*_2 \in \tilde{\mathcal{C}}_y $) is satisfied provided that
 $(y^*_1)'(x_0)> \max(1,\delta/r)$. This is automatically ensured by the assumptions
 made on $a,b,c$ since $(y^*_1)'(x_0)=(bx_0^2+2c
 x_0+ac)/(b^2x_0^2+2bcx_0+c^2)$.
That family is such that $\tsigma_1(y)=\tsigma_2(y)$ for $y \le y^*_2(x_0)$. 
 We can then calculate $\sigma_2$ as in subsection
 \ref{an_ex} using the relation $\sigma_2(x)=\frac{\sqrt{2(ry^*_2(x)- \delta
     x)(y^*_2(x)-x)}}{x \sqrt{y^*_2(x)'}}$. This gives
 $\sigma_2(x)=\sigma(x)$  for  $x \le x_0$  and for $x \ge x_0$,
\[
\sigma_2(x)=\sqrt{2\frac{[(r(y^*_1)'(x_0)-\delta ) x + r (y^*_1(x_0)-x_0(y^*_1)'(x_0))][((y^*_1)'(x_0)-1 ) x +  y^*_1(x_0)-x_0(y^*_1)'(x_0)]}{x^2 (y^*_1)'(x_0) }}. \]

\begin{figure}[h]
  \psfrag{cas1}{\tiny \hspace{-2.5cm} \boxed{x=0.5 \ $and$ \  y=0.4}}
  \psfrag{cas2}{\tiny \hspace{-2.5cm} \boxed{x=3 \ $and$ \ y=1} }
  
  \centerline{\psfig{file=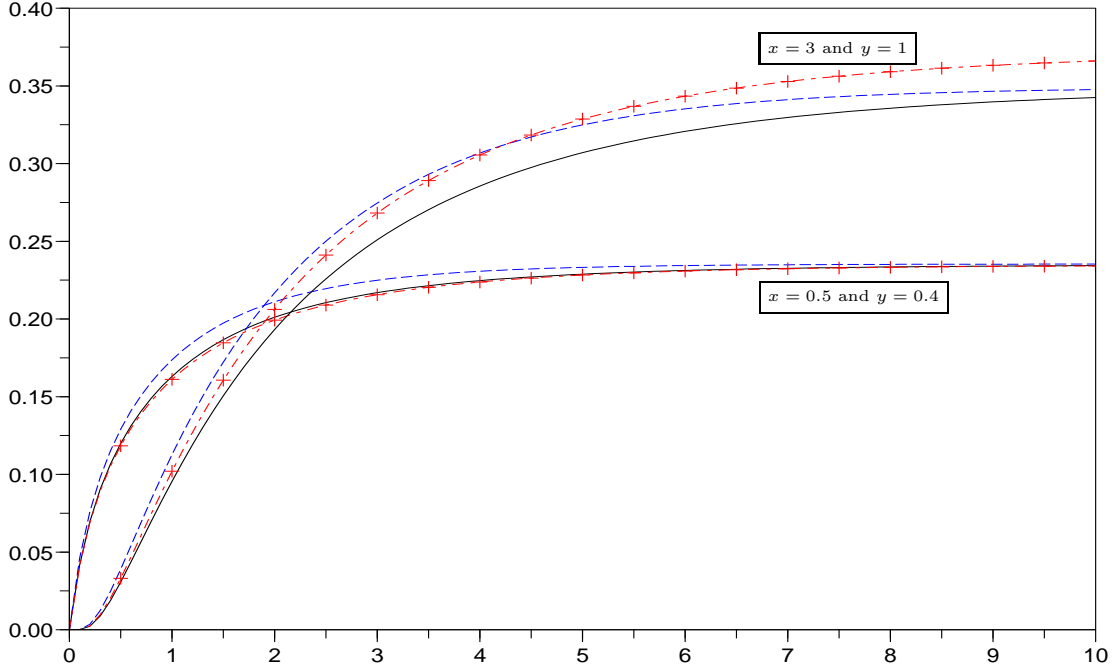,width=18cm,height=11cm} }
\caption{{\small $P_\sigma(T,x,y)$ (solid line), $P_{\sigma_2}(T,x,y)$ (dashed line
    with crosses) and $c_{\tsigma}(T,y,x)$ (dashed line) in function of the time
    $T$ for $a=1$, $b=0.4$, $c=0.1$, $x_0=0.5$, $r=0.2$ and $\delta=0.1$.}} \label{Dualfig2}
\end{figure} 

In Figure \ref{Dualfig2}, we have plotted the same example as in Figure
\ref{Dual_inf} ($x=0.5 \ $and$ \  y=0.4$), adding the graph of $T \mapsto P_{\sigma_2}(T,x,y)$. The volatility
function $\sigma_2$ has been calculated with the formula above with $x_0=0.5$.
According to Proposition \ref{equivalences} and the Duality, the three prices are
equal when $T$ is large. In the second example ($x=3 \ $and$ \  y=1$), we still
observe that $P_\sigma(T,x,y)$ and $c_{\tsigma}(T,y,x)$ converge toward the
same value when $T$ is large. On the contrary, the limit price of $P_{\sigma_2}(T,x,y)$ is significantly different. To observe the same price, we
should have taken, according to Proposition \ref{equivalences}, $x_0 \ge 3$.

\subsection{Calibration to ``real'' call and put prices}
In that subsection, we address some problems that arise if one tries to apply the
calibration procedure when the prices $p(K)$ and $c(K)$ do not derive from a
time-homogeneous model. We assume however that they are smooth functions of the
strike $K$, and focus for example on the calibration to put prices.

Firstly, let us observe that the arbitrage-free theory allows to define a dual volatility as previously by $(0,Y]$: \begin{equation}\label{defdual} \forall
K <Y, \
\eta_p(K)=\frac{1}{K} \sqrt{\frac{2(\delta p(K) +  K (r-\delta)
    p'(K))}{p''(K)}}.\end{equation}
Indeed, the payoff convexity in $K$ ensures the positivity of $p''(K)$ and  the arbitrage-free assumption ensures that $\delta p(K) +  K (r-\delta)
p(K)$ is nonnegative, so that the square-root is well defined. Let us prove the last
point and
 suppose the contrary (i.e. $\exists y>0$ such that $\frac{d}{dy}e^{\delta
  y}p(e^{(r-\delta)y})<0$) to exhibit an arbitrage opportunity. In that
case, there is $z>y$ such that $e^{\delta   y}p(e^{(r-\delta)y}) > e^{\delta
  z}p(e^{(r-\delta)z})$. We then sell one put with strike
$e^{(r-\delta)y}$ and buy $e^{\delta(z-y)}$ puts with strike $e^{(r-\delta)z}$. This initial
transaction generates a positive flow. The hedging works
as follows: naming $\tau$ the time at which the put sold is exercised, we have to pay
$e^{(r-\delta)y} - S_\tau$. In other words, we receive one share and
borrow $e^{(r-\delta)y}$ in cash. We keep this position until time $\tau+z-y$. At this time, we have exactly
$e^{\delta(z-y)}$ shares and puts with strike $e^{(r-\delta)z}$. Thus, we obtain at
least $e^{\delta(z-y)}e^{(r-\delta)z}=e^{(r-\delta)y}e^{r(z-y)}$ and we cancel the debt.
The next proposition gives sufficient conditions that allow to construct an
homogeneous volatility which is consistent to the observed prices.
\begin{proposition}\label{calib1put}
Let us assume that $K \in \mathbb{R}_+^* \mapsto p(K)$ is a
$\mathcal{C}^1$ function, $\mathcal{C}^2$ on
$\mathbb{R}_+^*-\{Y\}$ with $Y=\inf\{K>0:p(K)=K-x_0\}<+\infty$. Let us
also assume that $\eta_p$ defined by~(\ref{defdual}) is
bounded from below and above by two positive constants and admits a
left-hand limit in $Y$.
Then, if we extend $\eta_p$ in any continuous function on $(0,+\infty)$ satisfying
($\Hs$) still denoted by $\eta_p$, we have
\[\forall K>0, P_{\utilde{\eta_p}}(x_0,K)=p(K).\]
\end{proposition}
 Notice that once we choose the extended function $\eta_p$, we obtain
 $\utilde{\eta_p}$ by first solving~(\ref{EDy*}) on $\mathbb{R}_+^*$ starting
from $x_0$ with the condition $y^*_{\eta_p}(x_0)=Y$ and then using~(\ref{Dual_eta}).
\begin{proof}
The functions $K \mapsto p(K)$ and $K \mapsto c_{\eta_p}(K,x_0)$ solve
(\ref{EDg}). Since we have $0 \le p(K) \le K$ for arbitrage-free reasons, both
functions go to $0$ when $K \rightarrow 0$. Thanks to Remark \ref{LIMCOND}, they
are proportional to $g_{\uparrow}$ and therefore there is $\lambda >0$ such that:
\[ \forall K \le Y,\ p(K)=\lambda c_{\eta_p}(K,x_0). \]
The $\mathcal{C}^1$ assumption made on $p$ ensures $p(Y)=Y-x_0$ and $p'(Y)=1$. This
gives $g_{\uparrow}(Y)/g_{\uparrow}'(Y)=Y-x_0$ and therefore $Y=\ysd_{\eta_p}(x_0)$ using
Lemma \ref{charac_bound}. Thus,  $c_{\eta_p}(Y,x_0)=Y-x_0=p(Y)$ and
$\lambda=1$. One concludes with Theorem \ref{DualProp}.
\end{proof}

For the call case, everything  works in the same manner, but we need to assume moreover that
$c(K) \rightarrow 0$ when $K \rightarrow + \infty$. This is a rather natural
hypothesis that plays the same role as $p(K)\rightarrow 0$ when $K \rightarrow 0$.
\begin{proposition}\label{calib2put}
Let us assume that $K \in \mathbb{R}_+^* \mapsto c(K)$ is a
$\mathcal{C}^1$ function, $\mathcal{C}^2$ on
$\mathbb{R}_+^*-\{X\}$ with $X=\sup \{K>0, \ c(K)=x_0-K \}>0$ and $\lim_{K\rightarrow  + \infty}c(K)=0$.
Let us also assume that $\eta_c$ defined by \[ \forall K > X, \  \eta_c(K)=\frac{1}{K}
\sqrt{\frac{2 (\delta c(K) +  K (r-\delta)
    c'(K))}{c''(K)}}\]  is
bounded from below and above by two positive constants and admits a
right-hand limit in $X$.
Then, if we extend $\eta_c$ in any continuous function on $(0,+\infty)$ satisfying
($\Hs$) still denoted by $\eta_c$, we have
\[\forall K>0, C_{\uh{\eta_c}}(x_0,K)=c(K)\]
where $\uh{\eta_c}$ is obtained from $\eta_c$ like $\sigma$ from
$\hat{\sigma}$ in the end of the proof of Proposition \ref{calib}.
\end{proposition}

Therefore, we are able to find volatility functions that give exactly the put
prices and others that give exactly the call prices. Now, the natural question
is whether one can find a volatility function $\sigma$ 
that is consistent to both the put and call prices. According to
Proposition \ref{equivalences}, all the volatility functions
$\utilde{\eta_p}$ (resp. $\uh{\eta_c}$) giving the
put (resp. call) prices coincide on $(0,x_0)$
(resp. $(x_0,+\infty)$). The only volatility function
possibly giving both the put and call prices is
\begin{equation*}
   \sigma(x)=\begin{cases}
      \utilde{\eta_p}(x)\mbox{ if }x<x_0\\\uh{\eta_c}(x)\mbox{ if }x>x_0
   \end{cases}.
\end{equation*}
We deduce from Proposition \ref{equivalences} :
\begin{proposition}
  Assume that $\utilde{\eta_p}(x_0^-)=\uh{\eta_c}(x_0^+)$. Then, 
  \[\forall K>0,  p(K)=P_\sigma(x_0,K) \text{ and } c(K)=C_\sigma(x_0,K)\;\; \mbox{iff $\xss(Y)=x_0$ and $\ysp_\sigma(X)=x_0$.}\]
\end{proposition}

\section{The Black-Scholes model:  the unique model invariant through
this duality}

The purpose of that section is to put in evidence the particular role played by the
Black-Scholes' model for the perpetual American call-put duality. We have recalled in
the introduction that constant volatility functions are invariant by the duality. We
have also mentioned that for the European case, the call-put duality 
holds for all maturities without any change of the volatility function. Here, on the contrary,
we are going to prove that if the duality holds for the perpetual American options with the same volatility:
\begin{equation}\label{eq:BS_uni}
  \forall x,y>0 \ P_\sigma(x,y)=c_\sigma(y,x)  \end{equation}
then, under some technical assumptions, necessarily $\sigma(.)$ is a constant function.

\begin{proposition}
  Let us consider a positive interest rate $r$ and a nonnegative dividend rate
  $\delta<r$. We suppose that the volatility function $\sigma$ satisfies $(\Hs)$, and
  is analytic in a neighborhood of $0$, i.e.
  \begin{equation}\label{DSE}
    \exists \rho >0, \forall x \in [0,\rho ), \ \sigma(x)=\sum_{k=0}^{\infty}
    \sigma_k x^k.
  \end{equation}
  Then,~(\ref{eq:BS_uni}) holds if and only if $\forall x \ge 0, \ \sigma(x)=\sigma_0$.
\end{proposition}

We have already shown in the introduction that~(\ref{eq:BS_uni}) holds
in the Black-Scholes' case. So we only have to prove the necessary condition. We decompose the
proof into the three following lemmas.

\begin{lemma}\label{lemme_dual1}
Let us consider a volatility function that satisfies $(\Hs)$. If the dual volatility function $\tsigma$ is
 analytic in a neighborhood of $0$, then the boundaries $x^*_\sigma$ and $\ysd_{\tsigma}$ are also analytic in a neighborhood of $0$.
\end{lemma}
\begin{lemma}\label{lemme_dual2}
  Let us suppose that $\sigma$ satisfies $(\Hs)$ and is analytic in a neighborhood of $0$.
  Let us assume moreover that $r > \delta$.
  If the equality~(\ref{eq:BS_uni}) holds, $\sigma$ is constant in
  a neighborhood of $0$: 
  \[ \exists \rho >0, \forall y \in [0,\rho ], \ \sigma(y)=\sigma_0 .\]
\end{lemma}  
\begin{lemma}\label{lemme_dual3}
  Let us suppose that $\sigma$ is a constant function on $[0,\rho]$ for $\rho>0$
  satisfying $(\Hs)$ and~(\ref{eq:BS_uni}). Then, $\sigma$ is constant on
  $\mathbb{R}_+$ (and $x^*_\sigma$ and $\ysd_{\tsigma}$ are linear functions).
\end{lemma}

\begin{proof}[Proof of Lemma \ref{lemme_dual1}.] Let us first show that $\xss$ is analytic
in $0$. Thanks to the relation~(\ref{Bound_cond}), we have
$\frac{g(\ysd_{\tsigma}(x))}{g'(\ysd_{\tsigma}(x))}=\ysd_{\tsigma}(x)-x$, and therefore $\frac{g(y)}{g'(y)}=y-\xss(y)$.
Thus, $\xss(y)$ is analytic in $0$ iff $\phi(y)=\frac{g(y)}{g'(y)}$ is analytic in
$0$. Using the relation~(\ref{EDg}) and $\phi'=1-\frac{g''}{g'}\phi$, we get that
$\phi$ is solution of
\begin{equation}\label{EDO_phi}
  \phi'(y)=1+\frac{2}{\tsigma^2(y)} \left( (\delta-r) \phi(y)/y - \delta
  (\phi(y)/y)^2  \right).  \end{equation}
Notice that $\phi(y)=y-\xss(y)$ and \eqref{bound_boundaries_x} imply 
that if $\phi$ is analytic in $0$ then the coefficient of order $0$ in
  its expansion vanishes and the coefficient of order $1$ belongs to $(0,1)$.

To complete the proof we are first going to check that if $\psi(y)=\sum_{k=1}^{\infty} \phi_k
y^k$ with $\phi_1\in (0,1)$ solves \eqref{EDO_phi} in a neighborhood of
$0$ then $\phi\equiv \psi$ in this neighborhood. Then we will prove
existence of such an analytic solution $\psi$. We have $\psi(0)=0$, and
the function $\psi$ being analytic with $\phi_1\neq 0$, its
zeros are isolated points. There is therefore a 
neighborhood of $0$, $(0,2\epsilon)$ where $\psi$ does not vanish. Let us consider
$\gamma$ a solution of $\gamma'-\frac{1}{\psi} \gamma = 0$ starting from $\gamma(\epsilon) \not = 0$ in
$\epsilon$~: $\gamma(x)=\gamma(\epsilon) \exp \left( \int_{\epsilon}^x
  \frac{1}{\psi(u)}du \right) $. Since $\psi$ solves~(\ref{EDO_phi}), it is
not hard to check that $\gamma$ is solution of~(\ref{EDg}) with $\eta=\tsigma$. The limit condition $\gamma(x)
\underset{x \rightarrow 0}{\rightarrow} 0$ (cf. Remark \ref{LIMCOND}, still valid for
$g_\uparrow$ when $\delta=0$) is satisfied   
since we have $\frac{1}{\psi(u)} \underset{u \rightarrow 0} {\sim} \frac{1}{\phi_1
  u}$ and so $\int_{\epsilon}^x \frac{1}{\psi(u)}du \underset{x \rightarrow
  0}{\rightarrow} - \infty$.
Thus we have $\gamma(y)=c g(y)$ with $c \not = 0$ and $\psi(y)=g(y)/g'(y)=\phi(y)=y-\xss(y)$. We can
then write $\xss(y)= (1-\phi_1)y-\sum_{k=2}^{\infty} \phi_k y^k$ in the
neighborhood of $0$
with $1-\phi_1>0$. It is well-known that in that case, the reciprocal
function $\ysd_{\tsigma}$ is also analytic in $0$.

Let us turn to the existence of $\psi$. Since $\sigma_0 \ge
\underline{\sigma}>0$, $y\rightarrow \frac{2}{\tsigma^2(y)}$ is an analytic
function in the neighborhood of $0$. Thus, there is $\rho_0>0$ and $a_0>0$ such that
\[\forall y \in [0,\rho_0], \  \frac{2}{\tsigma^2(y)}=\sum_{k=0}^{\infty} a_k y^k
\text{ and } \sum_{k=0}^{\infty} |a_k| \rho_0^k < \infty
. \]
The analytic function $\sum_{k\geq 1}\phi_k y^k$ solves~(\ref{EDO_phi}) if and only if
 \[ \sum_{k=0}^{\infty} (k+1) \phi_{k+1} y^k = 1 + (\delta-r) \sum_{k=0}^{\infty}
\left(\sum_{i+j=k} a_i \phi_{j+1} \right) y^k - \delta \sum_{k=0}^{\infty}
\left(\sum_{i+j+l=k} a_i \phi_{j+1}\phi_{l+1} \right) y^k.\]
Identifying the terms of order $0$, we get that $\phi_1$ solves
$P(\phi_1)=0$ where $P(x)=\delta a_0 x^2+(1-(\delta-r)a_0)x-1$. Since
$P(0)=-1<0$ and $P(1)=ra_0>0$, the polynomial $P$ admits a unique root on
$(0,1)$ and we choose $\phi_1$ equal to this root. Then, by
identification of
the terms with order $k$, we define the sequence $(\phi_k)_{k\geq 1}$
inductively by
\[ \phi_{k+1} = \frac{ (\delta-r)\sum_{i+j=k,j \not = k} a_i \phi_{j+1} -
\delta \sum_{i+j+l=k, j \not = k, l \not = k} a_i \phi_{j+1}\phi_{l+1}  }{k+1
+ (r-\delta)a_0+2 \delta a_0 \phi_1}. \]
This ratio is well defined since $(r-\delta)a_0+2 \delta a_0 \phi_1 =\delta a_0
\phi_1 + 1/\phi_1- 1 >0$.\\
We still have to check that the series $\sum_{k\geq 1}\phi_k y^k$ is
defined in a neighborhood of $0$. To do so, we are going to show that
there is $\rho >0 $ such that the sequence 
$(|\phi_k| \rho^k)_{k\geq 1}$ is bounded. We have for $1 \le k \le n$:
\[|\phi_{k+1}| \rho^{k}  \le \frac{ |\delta-r| \underset{j=0}{\overset{k-1}{\sum}} |a_{k-j}| \rho^{k-j}
  |\phi_{j+1}|\rho^{j} + \delta \underset{i=0}{\overset{k}{\sum}}
  \left( \underset{j \not = k, l \not = k }{\underset{j+l=k-i}{\sum}}   |\phi_{j+1}|\rho^{j}
    |\phi_{l+1}|\rho^{l} \right) |a_{i}| \rho^{i} }{k+1}. \]
Let us suppose that for $1 \le j <k$,  $|\phi_{j+1}|\rho^{j} \le 1/(j+1)$. Then,
\[|\phi_{k+1}| \rho^{k}  \le\frac{|\delta-r| \rho \sum_{j=1}^k |a_j| \rho^{j-1} +
  \delta \underset{i=0}{\overset{k}{\sum}}  \left( \underset{j+l=k-i}{\sum}
    \frac{1}{j+1} \frac{1}{l+1}\right)|a_{i}| \rho^{i} }{k+1}  . \]
We remark that $\underset{j+l=k-i}{\sum}
\frac{1}{j+1} \frac{1}{l+1}=\frac{1}{k-i+2}\underset{j+l=k-i}{\sum} \frac{1}{j+1}
+ \frac{1}{l+1} \le 2 \frac{\ln(k-i+1)+1}{k-i+2}$, and we finally get:
\begin{equation}\label{lemme1::1}
  |\phi_{k+1}| \rho^{k}  \le \frac{2 \delta |a_0|\frac{\ln(k+1)+1}{k+2} + \rho (|\delta-r|+2\delta) \sum_{j=1}^k |a_j| \rho^{j-1} }{k+1}
\end{equation}
since $\frac{\ln(k-i+1)+1}{k-i+2} \le 1$. Let us now consider $k_0$ such that $\forall
k \ge k_0, 2 \delta |a_0|\frac{\ln(k+1)+1}{k+2}<1/2$. Now, we chose $\rho\in (0,\rho_0)$ small enough such
that $\forall k\le k_0, |\phi_{k+1}|\rho^{k} \le 1/(k+1)$ and $\rho
(|\delta-r|+2\delta) \sum_{j=1}^{\infty} |a_j| \rho^{j-1} < 1/2$. Then we get by
induction from~(\ref{lemme1::1}) that $\forall k \ge k_0, |\phi_{k+1}|\rho^{k} \le
1/(k+1)$.\end{proof}

\begin{proof}[Proof of Lemma \ref{lemme_dual2}.]
  On the one hand, thanks to the assumption, $\sigma=\tsigma$ is analytic in 0, and
  therefore $\xss$ is analytic in 0 thanks to Lemma \ref{lemme_dual1}:
 \[ \exists \rho >0, \ \forall y \in [0,\rho), \xss(y)=\sum_{i=1}^{\infty} 
x_i y^i \text{ and } \sigma(y)=\sum_{i=0}^{\infty} \sigma_i y^i.\]
On the other hand, it is  not hard then to deduce from~(\ref{Dual_sig}),  $\sigma=\tsigma$ and the 
differential equation~(\ref{EDx*}) that
\begin{equation}\label{EDx*ter}
  \xss(y)'=\frac{2(y-\xss(y))(ry-\delta \xss(y))  }{y^2 \sigma(y)^2}.
\end{equation}
{From} Corollary \ref{sig->x*} and~(\ref{sigma_de_x*}), we get 
\begin{equation}\label{EDx*dual} \xss(y)'=\frac{(y-\xss(y))(ry-\delta \xss(y)) ((\xss)^{-1})'(y) }{ ( (\xss)^{-1}(y) - y) (
  r(\xss)^{-1}(y) -\delta y) } .\end{equation} 
Now, we consider $n =\inf \{i \ge 2, x_i \not =0 \}$ and suppose it finite. We can
get easily that:

\begin{tabular}{ll}
  $\xss(y)=x_1 y +x_n y^n+\dots$ & $\xss(y)'=x_1  + n x_n y^{n-1}+\dots$ \\
  $(\xss)^{-1}(y)=\frac{1}{x_1} y - \frac{x_n}{x_1^{n+1}} y^n+\dots$ &
  $((\xss)^{-1})'(y)=\frac{1}{x_1}(1  - \frac{nx_n}{x_1^{n}} y^{n-1})+\dots$ \\
\end{tabular}\\
and then
\begin{eqnarray*}
(1-\xss(y)/y) (r-\delta \xss(y)/y) & = & (1-x_1)(r-\delta x_1) +x_n(2\delta x_1 -(r+\delta))y^{n-1}+\dots \\
\bigg( \frac{(\xss)^{-1}(y)}{y} - 1\bigg) \bigg(  r\frac{(\xss)^{-1}(y)}{y} -\delta \bigg) & = & \frac{1}{x_1^2}\left\{
  (1-x_1)(r-\delta x_1) +\frac{x_n}{x_1^{n}} ((r+\delta)x_1-2r)y^{n-1} \right\} +
\dots 
\end{eqnarray*}
The right hand side of~(\ref{EDx*dual}) has then the following expansion:
\[ x_1\left\{ 1+ \frac{x_n}{(1-x_1)(r-\delta x_1)}
  \left[ 2 \delta x_1 -(r+\delta)
    +\frac{2r}{x_1^n}-\frac{r+\delta}{x_1^{n-1}}\right]y^{n-1} - \frac{nx_n}{x_1^{n}}y^{n-1}\right\}+
\dots 
\]
The equality of the terms of order $n-1$ in~(\ref{EDx*dual}) then leads to:
\[nx_n x_1^{n-1}=\frac{x_n}{(1-x_1)(r-\delta x_1)} \left[ 2 \delta x_1^{n+1} -(r+\delta)x_1^{n}
-(r+\delta)x_1 +2r \right] -nx_n. \]
Since $x_n \not = 0$ and with a simplification we get
\begin{equation}\label{eqx1}
  n (1+x_1^{n-1})=\frac{1}{r-\delta x_1}\left[-2 \delta x_1^n +(r-\delta)\sum_{k=1}^{n-1}x_1^{k}+2r\right].
\end{equation}
In the case $\delta=0$ this gives $n (1+x_1^{n-1})=x_1^{n-1}+\dots+x_1+2$ which is
not possible because $x_1 \in (0,1)$. When $0< \delta <r$, we denote $\alpha=r/
\delta >1$ and rewrite~(\ref{eqx1}):
\begin{equation*}
  n (1+x_1^{n-1})(\alpha - x_1)=-2  x_1^n
  +(\alpha-1)x_1^{n-1}+\dots+(\alpha-1)x_1+2\alpha=\alpha-x_1^n +
  (\alpha-x_1)\frac{1-x_1^n}{1-x_1}.
\end{equation*}
Therefore, $ n (1+x_1^{n-1})=\frac{\alpha-x_1^n}{\alpha-x_1}+\frac{1-x_1^n}{1-x_1}
<2\frac{1-x_1^n}{1-x_1}$ because $\beta \mapsto
\frac{\beta-x_1^n}{\beta-x_1}$ is decreasing on $(1,\alpha)$ ($x_1^n<x_1$). To show that this
is impossible, we consider $P_n(x)= n (1+x^{n-1})-2\sum_{k=0}^{n-1} x^k$. We have
$P_n(1)=0$ and for $x < 1$,  $P'_n(x)=n(n-1)x^{n-2}-2 \sum_{k=1}^{n-1}k x^{k-1}=2
\sum_{k=1}^{n-1}k (x^{n-2}-x^{k-1}) < 0$. Thus $P_n$ is positive on $[0,1)$ and
$P_n(x_1)>0$ which is a contradiction. \end{proof}

\begin{proof}[Proof of Lemma \ref{lemme_dual3}.]
  It is easy to get from~(\ref{EDx*}) and $\sigma=\tsigma$ that
\begin{equation}\label{EDx*bis}
     \xss(y)'=\frac{ \xss(y) \sigma( \xss(y) ) }{y \sigma(y)}.
\end{equation}
We have $\sigma(x)=\sigma_0$ for $x \in [0,\rho]$. Since $\xss(y)$  solves~(\ref{EDx*bis}) and
$\xss(y) \le y$, $\xss(y)'=\xss(y) /y$ on $[0,\rho]$. Therefore, $\xss(y)=x_1 y$ for 
$y \in [0,\rho]$. Thanks to~(\ref{EDx*}), $x_1$ is the unique root in $(0,\min(1,r/\delta))$ of
\[ x_1 \sigma_0^2= 2(1-x_1)(r-\delta x_1 ). \]
Now let us observe that~(\ref{EDx*}) gives for $y \in (0,\ysd_{\tsigma}(\rho)]$, 
$\xss(y)'= \frac{\xss(y)^2 \sigma_0^2}{2(y-\xss(y))(ry-\delta \xss(y))}$
with $\xss(\rho)=x_1\rho$. Since $y\rightarrow x_1y$ solves this ODE,
for which local uniqueness holds thanks to the Cauchy Lipschitz theorem,
we then have
$\xss(y)=x_1 y$ on $[\rho,\ysd_{\tsigma}(\rho)]$ and so $\ysd_{\tsigma}(\rho)=(\xss)^{-1}(\rho)=\rho/x_1$.
Then,~(\ref{EDx*bis}) gives $\sigma_0/\sigma(y)=1$ on
$[\rho,\rho/x_1]$. Thus, we prove by induction on $n$ that $\xss(y)=x_1 y$ and
$\sigma(y)=\sigma_0$ for $y \in [0,\rho/(x_1)^n]$. This shows the desired result.
\end{proof}

\section{Conclusions and further developments}

Addressing Call-Put duality for American options with finite maturity in
models with time-dependent local volatility functions like
\eqref{intedsprim} would be of great interest. For the perpetual case
treated in this paper, we could take advantage of a very nice feature :
in the continuation region, the price of the option writes as the
product of a function of the underlying spot price by another function
of the strike price. Unfortunately, this product property no longer
holds in the general case.

Next, according to our numerical experiments
(see figure \ref{Dualfig2}), American Put and Call prices computed in infinite maturity
dual models may differ for finite maturities.
This means that in the case of a time-homogeneous primal local volatility function
$\varsigma(t,x)=\sigma(x)$, if there exists a dual local volatility function
for some finite maturity $T$, then this volatility function is either
time-dependent or depends on the maturity $T$. On the contrary, in the European case
presented in the introduction, time-homogeneous volatility functions are
preserved by the duality.

Let us nevertheless conclude on an encouraging remark. Let $P(T,x,y)$
denote the initial price of the American Put option with maturity $T$
and strike $y$ in the model \eqref{intedsprim} and $x^*(T,y)$ stand for
the corresponding exercise boundary such that $P(T,x,y)=(y-x)^+$ if and
only if $x\leq x^*(t,y)$. Then the smooth-fit principle writes
\begin{equation*}
  \begin{cases}
     P(T,x^*(T,y),y)=y-x^*(T,y)\\
\partial_x P(T,x^*(T,y),y)=-1
  \end{cases}.
\end{equation*}
Differentiating the former equality with respect to $y$ yields
$$\partial_xP(T,x^*(T,y),y)\partial_yx^*(T,y)+\partial_yP(T,x^*(T,y),y)=1-\partial_yx^*(T,y).$$
With the second equality, one deduces that
$\partial_yP(T,x^*(T,y),y)=1$. Therefore the smooth-fit principle
automatically holds
for the dual Call option if there exists any. 


\end{document}